\documentclass [10pt,reqno]{amsart}

\usepackage {hannah}

\usepackage{color}
\usepackage{graphicx}
\usepackage{texdraw}
\usepackage{epic}
\usepackage{amsmath}
\usepackage{amssymb} 
\usepackage{latexsym}
\usepackage{bbm}

\DeclareMathOperator{\val}{val}
\DeclareMathOperator{\ini}{in}
\DeclareMathOperator{\trop}{trop}
\DeclareMathOperator{\Trop}{Trop}
\DeclareMathOperator{\MaxTrop}{MaxTrop}
\DeclareMathOperator{\Int}{Int}
\DeclareMathOperator{\conv}{conv}
\DeclareMathOperator{\tdim}{tdim}
\newcommand{\bx}{{\boldsymbol x}}
\newcommand{\bp}{{\boldsymbol p}}
\newcommand{\bm}{m}

\DeclareMathAlphabet{\mathds}{U}{dsrom}{m}{n}
\newcommand{\PP}{{\mathds P}}

\title [Tropical curves with a singularity in a fixed point]{Tropical curves with a singularity in a fixed point}

\author {Hannah Markwig}
\address {Hannah Markwig,  Fachrichtung Mathematik, 
   Universit\"at des Saarlandes, 66041 Saar\-br\"ucken, Germany}
\email {hannah@math.uni-sb.de}

\author{Thomas Markwig}
\address{Thomas Markwig, Fachbereich Mathematik, TU Kaiserslau\-tern, 
   67653 Kaiserslautern, Germany}
\email {keilen@mathematik.uni-kl.de}

\author{Eugenii Shustin}
\address{Eugenii Shustin, School of Mathematical Sciences, Tel Aviv University, Ramat Aviv, Tel Aviv 69978, Israel}
\email {shustin@post.tau.ac.il}

\thanks{The authors were supported by the
  Hermann-Minkowski Minerva Center for Geometry at the Tel Aviv
  University and by the Mathematical Sciences Research Institute
  (MSRI).
  The first and second author were also supported by the German Research Foundation
  (Deutsche Forschungsgemeinschaft (DFG)) through
  the Institutional Strategy of the University of G\"ottingen.
  The third author also enjoyed a support from the Israeli
  Science Foundation grant no. 448/09.}

\begin {document}

   \maketitle

   \section{Introduction}
The aim of this paper is to start the investigation of tropical singularity theory, i.e.\ the study of the geometry of singular tropical varieties and their families. Singular plane tropical curves (notably, nodal plane tropical curves) appear as the main technical tool in the tropical computation of the Gromov-Witten invariants of $\mathbb{P}^2$ (\cite{Mi03}). However, no systematic theory of tropical singularities has been built, and our work is a first step towards such a theory. Namely, we consider the simplest possible case: plane tropical curves with a given Newton polygon $\Delta$ which are tropicalisations of algebraic curves having a singularity at a given point. 
That is, we study the tropicalisations of algebraic curves on the toric surface corresponding to $\Delta$ with a singularity at a given point in the torus. 
We construct the tropical linear space $\Trop(\Sing_\K(\Delta))$ which is the tropicalisation of 
the family of plane curves with a given Newton polygon $\Delta$ and a singularity at a given point and we classify 
its maximal cones. We show that the parameterising space $\MaxTrop(\Sing_\K(\Delta))$ for the tropical curves under consideration is a connected 
in codimension one proper polyhedral subcomplex of $\Trop(\Sing_\K(\Delta))$ of the same dimension. 
We classify the tropical curves which are generic elements of the top-dimensional cones of $\MaxTrop(\Sing_\K(\Delta))$; furthermore, lifting such tropical curves into algebraic ones, we describe possible 
refined tropical limits and construct canonical simple parameterisations of the considered tropical curves.

As Newton polygon, fix a non-degenerate convex
   lattice polygon  $\Delta\subset\R^2$ and denote by
   $\mathcal{A}=\Delta\cap\Z^2$ the lattice points of $\Delta$. For
   any field $\K$ we denote by 
   $\Sigma=\Tor_\K(\Delta)$ the toric surface associated to $\Delta$. It comes with
   the tautological line bundle ${\mathcal L}_\Delta$
   generated by the global sections $\{x^iy^j\ :\
   (i,j)\in\mathcal{A}\}$. The torus $(\K^*)^2$ is embedded in
   $\Sigma$ via
   \begin{displaymath}
     \Psi_{\mathcal{A}}:(\K^*)^2\longrightarrow \PP_\K^{\mathcal{A}}:(x,y)\mapsto
     \big(x^iy^j\;|\;(i,j)\in\mathcal{A}\big)
   \end{displaymath}
   and inside the torus the elements in the linear system $|{\mathcal
     L}_\Delta|$ are defined by the equations
   \begin{displaymath}
     f_{\underline{a}}=\sum_{(i,j)\in\mathcal{A}}a_{i,j}\cdot x^i\cdot
     y^j=0.
   \end{displaymath}
   $|{\mathcal L}_\Delta|$
   contains a nonempty linear subsystem $\Sing_\K(\Delta)$ of curves with
   a singularity at the point $\bp=(1,1)$. This
   subsystem is given by the linear equations
   \begin{displaymath}
     f_{\underline{a}}(\bp)=0,\;\; \frac{\partial f_{\underline{a}}}{\partial
       x}(\bp)=0,\;\;\frac{\partial f_{\underline{a}}}{\partial
       y}(\bp)=0.
   \end{displaymath} More details about this family follow in Section \ref{sec-family}.
   Our first aim is to describe
   the geometry of $\Trop(\Sing_\K(\Delta))$, the tropicalisation of
   $\Sing_\K(\Delta)$, as a tropical variety (i.e.\ a balanced fan) in
   $\R^{s-1}=\R^{\mathcal{A}}/(1,\ldots,1)\cdot\R$ (where
   $s=\#\mathcal{A}$) and to analyse the underlying tropical curves
   (i.e.\ the tropicalisations of the singular algebraic curves $C\in\Sing_\K(\Delta)$).

   In order to be able to tropicalise we use an algebraically
   closed field $\K$ with a valuation  $\val:\K^*\longrightarrow\R$  whose value group is dense in $\R$, e.g.\ the field of Puiseux series.
  For an ideal $I\subset \K[x_1,\ldots,x_n]$ determining an affine variety
   $V=V(I)\subset \K^{n}$ we define the \emph{tropicalisation} of $V$ to be
   \begin{displaymath}
     \Trop(V):= \overline{\{(-\val(x_1),\ldots,-\val(x_n)) \;|\;
       (x_1,\ldots,x_n)\in V(I) \cap (\K^{\ast})^n\}}.
   \end{displaymath}
For an introduction to tropicalisations of algebraic varieties, see e.g.\ \cite{RST03} and \cite{MS09}.
   If the ideal $I$ is homogeneous and defines a projective variety,
   we may consider $\Trop(V)$ modulo the linear space
   $(1,\ldots,1)\cdot\R$, i.e.\ we
   identify $\Trop(V)$ with its image in $\R^n / (1,\ldots,1)\cdot \R$.
If $I=\langle f
     \rangle \subset \K[x,y]$ and $f=\sum a_{ij}
       x^iy^j$, then the \emph{plane tropical curve} $\Trop(V(f))$ equals the locus of
     non-differentiability of the \emph{tropical polynomial}
 $ \trop f:= \max \{ -\val(a_{ij})+ix+jy\}
     $ by Kapranov's Theorem (see
     \cite[Theorem~2.1.1]{EKL06}). 
Plane tropical curves are dual to marked subdivisions of the Newton polygon of $f$  (see e.g.\ \cite[Prop.\ 3.11]{Mi03}), and their combinatorics can be studied in terms of these subdivisions. More details follow in the Section \ref{sec-secfan}, which introduces marked subdivisions, the secondary fan and its connection to plane tropical curves.

 Let us now give an example for the tropicalisation of a singular
   curve.

   \begin{example}
     We consider the polynomial
     \begin{displaymath}
       f=xy^2-tx^2-(2+t^3)\cdot xy+(1+2t+t^3)\cdot x+t^3y-(t+t^3)\in\K[x,y],
     \end{displaymath}
     where $\K$ is the field of Puiseux series.
     One easily verifies that $\bp=(1,1)\in(\K^*)^2$ is a singular
     point of the curve $V(f)$. $f$ defines a curve in the toric
     surface $\Tor_\K(\Delta)$, where $\Delta$ is as in the left hand
     side of Figure  \ref{fig:exsing}.
     \begin{figure}[h]
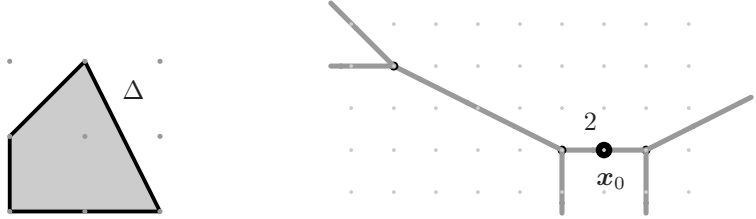

       \centering
       \begin{texdraw}
         \drawdim cm  \relunitscale 1
         \linewd 0.05
         \move (1 2)
         \lvec (2 0)
         \lvec (0 0)
         \lvec (0 1)
         \lvec (1 2)
         \lfill f:0.8
         \move (0 0) \fcir f:0.6 r:0.03
         \move (0 1) \fcir f:0.6 r:0.03
         \move (0 2) \fcir f:0.6 r:0.03
         \move (1 0) \fcir f:0.6 r:0.03
         \move (1 1) \fcir f:0.6 r:0.03
         \move (1 2) \fcir f:0.6 r:0.03
         \move (2 0) \fcir f:0.6 r:0.03
         \move (2 1) \fcir f:0.6 r:0.03
         \move (2 2) \fcir f:0.6 r:0.03
         \htext (1.5 1.5){$\Delta$}
       \end{texdraw}
       \hspace{2cm}
       \begin{texdraw}
         \drawdim cm  \relunitscale 0.7 \arrowheadtype t:V
         \linewd 0.1  \lpatt (1 0)
         \setgray 0.6
         \relunitscale 0.8
         \move (3 -1) \fcir f:0 r:0.1
         \move (3 -1) \lvec (1 -1)
         \htext (1.5 -0.5){$2$}
         \move (3 -1) \rlvec (2.5 1.25)
         \move (3 -1) \rlvec (0 -1.5)
         \move (1 -1) \fcir f:0 r:0.1
         \move (1 -1) \lvec (-3 1)
         \move (1 -1) \rlvec (0 -1.5)
         \move (-3 1) \fcir f:0 r:0.1
         \move (-3 1) \rlvec (-1.5 0)
         \move (-3 1) \rlvec (-1.5 1.5)
         \move (2 -1) \fcir f:0 r:0.2
         \htext (1.8 -1.9) {$\bx_0$}
         \move (-4 -2) \fcir f:0.8 r:0.05
         \move (-4 -1) \fcir f:0.8 r:0.05
         \move (-4 0) \fcir f:0.8 r:0.05
         \move (-4 1) \fcir f:0.8 r:0.05
         \move (-4 2) \fcir f:0.8 r:0.05
         \move (-3 -2) \fcir f:0.8 r:0.05
         \move (-3 -1) \fcir f:0.8 r:0.05
         \move (-3 0) \fcir f:0.8 r:0.05
         \move (-3 1) \fcir f:0.8 r:0.05
         \move (-3 2) \fcir f:0.8 r:0.05
         \move (-2 -2) \fcir f:0.8 r:0.05
         \move (-2 -1) \fcir f:0.8 r:0.05
         \move (-2 0) \fcir f:0.8 r:0.05
         \move (-2 1) \fcir f:0.8 r:0.05
         \move (-2 2) \fcir f:0.8 r:0.05
         \move (-1 -2) \fcir f:0.8 r:0.05
         \move (-1 -1) \fcir f:0.8 r:0.05
         \move (-1 0) \fcir f:0.8 r:0.05
         \move (-1 1) \fcir f:0.8 r:0.05
         \move (-1 2) \fcir f:0.8 r:0.05
         \move (0 -2) \fcir f:0.8 r:0.05
         \move (0 -1) \fcir f:0.8 r:0.05
         \move (0 0) \fcir f:0.8 r:0.05
         \move (0 1) \fcir f:0.8 r:0.05
         \move (0 2) \fcir f:0.8 r:0.05
         \move (1 -2) \fcir f:0.8 r:0.05
         \move (1 -1) \fcir f:0.8 r:0.05
         \move (1 0) \fcir f:0.8 r:0.05
         \move (1 1) \fcir f:0.8 r:0.05
         \move (1 2) \fcir f:0.8 r:0.05
         \move (2 -2) \fcir f:0.8 r:0.05
         \move (2 -1) \fcir f:0.8 r:0.05
         \move (2 0) \fcir f:0.8 r:0.05
         \move (2 1) \fcir f:0.8 r:0.05
         \move (2 2) \fcir f:0.8 r:0.05
         \move (3 -2) \fcir f:0.8 r:0.05
         \move (3 -1) \fcir f:0.8 r:0.05
         \move (3 0) \fcir f:0.8 r:0.05
         \move (3 1) \fcir f:0.8 r:0.05
         \move (3 2) \fcir f:0.8 r:0.05
         \move (4 -2) \fcir f:0.8 r:0.05
         \move (4 -1) \fcir f:0.8 r:0.05
         \move (4 0) \fcir f:0.8 r:0.05
         \move (4 1) \fcir f:0.8 r:0.05
         \move (4 2) \fcir f:0.8 r:0.05
       \end{texdraw}
       \caption{Tropicalisation of a singular curve}
       \label{fig:exsing}
     \end{figure}
     The tropicalisation of $V(f)$ is shown in the right hand side of
     Figure \ref{fig:exsing}. The singularity $\bp=(1,1)$ tropicalises
     to $\bx_0=(0,0)$. It sits precisely in the middle of an edge of
     weight two, i.e.\ it has the same distance to both neighbouring
     vertices. This metric condition and the fact that the edge has
     weight two are no coincidences, they are a general phenomenon, as
     we will see in Section \ref{sec-maxclass}.
   \end{example}

   Since $\Sing_\K(\Delta)$ is given by a linear ideal (see also Section
   \ref{sec-family}), we can use results of \cite{AK06} and
   \cite{FS05} to study its tropicalisation
   $\Trop(\Sing_\K(\Delta))$. We repeat
   facts about tropicalisations of linear ideals in Subsection
   \ref{subsec-Bergman} and study the top-dimensional cones of
   $\Trop(\Sing_\K(\Delta))$ in \ref{subsec-class}. We relate these
   top-dimensional cones to cones of the secondary fan in Subsection
   \ref{subsec-familyandsecfan}.
   We study the connection of $\Trop(\Sing_\K(\Delta))$ to the
   tropical discriminant in \ref{subsec-discriminant}.

Since many tropical polynomials can induce the same tropical curve, the family $\Trop(\Sing_\K(\Delta))$ is not the parameter space of singular tropical curves. Given a combinatorial type of a tropical curve, the set of all tropical curves of this type is naturally parametrised by a polyhedron whose dimension we call the type-dimension, $\tdim$.
We call a marked subdivision \emph{of maximal dimensional type} if the dimension of its corresponding cone in the secondary fan equals the $\tdim$ of the corresponding type (see Subsection \ref{subsec-dimofcurves}).
In Section \ref{sec-maxclass}, we study the subcomplex $\MaxTrop(\Sing_\K(\Delta))$ of $\Trop(\Sing_\K(\Delta))$ corresponding to maximal dimensional types. 
This subcomplex is the parameter space of singular tropical curves.
We classify the top-dimensional cones of $\MaxTrop(\Sing_\K(\Delta))$. A top-dimensional cone corresponds to a weight class of the lattice of flats of the matroid associated to the linear subspace $\Sing_\K(\Delta)$. We relate such a weight class with a specific neighbourhood of the singular point of the corresponding tropical curves. As the main consequence of this classification, we  show that the singularity of a tropical curve of maximal
   dimensional type is either a crossing of two edges of weight one, or a $3$-valent
   vertex of multiplicity $3$ adjacent to three edges of weight one, or a point on an edge of weight $2$ whose
   distances to the neighbouring vertices satisfy a certain metric
   condition.
A vertex of multiplicity $3$ adjacent to three edges of weight one is dual to a polygon with one interior point. This polygon could be further subdivided such that the dual picture is a cycle. The vertex can thus be viewed as the limit of a vanishing cycle, and we call it a \emph{nodal vertex}.

 The result of the classification, which is the main result of our paper, can be found in Theorem \ref{thm-class}.

Note that our result does not depend (up to shift) on the choice of the singular
   point, as long as it is a point in the torus $(\K^\ast)^2$ (see Remark \ref{rem-whichpoint} and \ref{rem-whichpoint2}).
We study what happens
   if we move the point to a coordinate line in Section
   \ref{sec-nottorus}.

In Section \ref{sec-alglift}, we study algebraic preimages, i.e.\ lifts, of our singular tropical curves,
   and in particular we thus give  a
   conceptual explanation for the metric condition mentioned
   above. The lift of a singular tropical curve of maximal dimensional type enhances the tropical curve with a refined tropical limit, which in particular allows us to define a canonical  $3$-valent parametrisation of the given tropical curve such that the tropical singular point lifts either to a nodal $3$-valent vertex, or to a pair of interior points of two edges, or to one or two $1$-valent vertices. This yields a more refined classification of singular tropical curves.

   We would like to thank Eva Maria Feichtner, Anders Jensen, Joaquim Ro\'{e}, Tim R\"omer and Kirsten
   Schmitz for valuable discussions, and we would like to thank the
   referee for his helpful remarks towards the presentation of the results.

   \section{The secondary fan and its relation to plane tropical curves}
   \label{sec-secfan}

   Here, we briefly repeat some basic definitions. For more details,
   see \cite[Chapter 7]{GKZ} or \cite{Mi03}.
   In the case of plane tropical curves, we can conclude from
   Kapranov's theorem that $\Trop(V(f))$ is a piece-wise linear
   graph in $\R^2$. An important fact is that this graph is \emph{dual} to a subdivision
   of the \emph{Newton polygon} $\Delta=\conv\{(i,j) \;|\; a_{ij}\neq 0\}$ of $f$.

   A \emph{marked polygon} $(Q,\mathcal{B})$ is a $2$-dimensional convex lattice
   polygon $Q$ in $\R^2$ together with a
   subset $\mathcal{B}$ of the lattice points $Q\cap\Z^2$ containing the
   vertices of $Q$.

   A \emph{marked subdivision} of a polygon
   $\Delta$ is a collection of marked polygons,
   $T=\{(Q_1,\mathcal{A}_1),\ldots,(Q_k,\mathcal{A}_k)\}$,
   such that
   \begin{itemize}
   \item $\Delta=\bigcup_{i=1}^k Q_i$,
   \item $Q_i\cap Q_j$ is a face (possibly empty) of $Q_i$ and of $Q_j$
     for all $i,j=1,\ldots,k$, and
   \item $\mathcal{A}_i\cap(Q_i\cap Q_j)=\mathcal{A}_j\cap (Q_i\cap Q_j)$ for all $i,j=1,\ldots,k$.
   \end{itemize}
   We do not require that $\bigcup_{i=1}^k \mathcal{A}_i=\Delta\cap \Z^2$.

   \begin{definition}
     We define the \emph{type} of a marked subdivision to be the
     subdivision, i.e.\ the collection of the $Q_i$, without the
     markings.
   \end{definition}

\begin{figure}[h]
     \centering
     \input{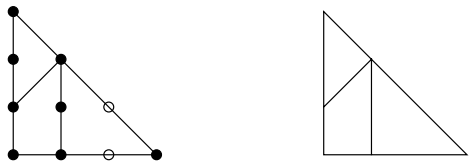}
     \caption{A marked subdivision and its type}
     \label{fig:type}
   \end{figure}

   Figure \ref{fig:type} shows an example of a marked subdivision and its type.
   The subset of lattice points which are marked in each
   $Q_i$ are drawn in black, the lattice points $\Delta\cap \Z^2$
   which are not marked are white. We will stick to this convention
   throughout the paper.

   For a finite subset $\mathcal{A}$ of the lattice $\Z^2$ we denote
   by $\R^{\mathcal{A}}$ the set of vectors indexed by the lattice
   points in $\mathcal{A}$.
   A point $u\in\R^{\mathcal{A}}$ induces a \emph{marked subdivision}
   of $\Delta$ by considering the convex hull of
   \begin{equation}\label{e1}
     \big\{(i,j,u_{ij})\;\big|\;(i,j)\in\mathcal{A}\}\subset\R^3
   \end{equation}
   in $\R^3$, and projecting the upper faces onto the $xy$-plane. A
   lattice point $(i,j)$ is marked if the point $(i,j,u_{ij})$ is
   contained in one of the upper faces. Marked subdivisions of
   $\Delta$ obtained in this way are called \emph{regular} or
   \emph{coherent}. We say two points $u$ and $u'$ in $\R^{\mathcal{A}}$
   are equivalent if and only if they induce the same regular marked subdivision
   of $\Delta$.  This defines an equivalence relation on
   $\R^{\mathcal{A}}$ whose equivalence
   classes are the relative interiors of convex cones. The collection of these cones is
   the \emph{secondary fan} of $\Delta$.

   Regular marked subdivisions of $\Delta$ are dual to plane tropical curves
   (see e.g.\ \cite[Prop.\ 3.11]{Mi03}).
   Given a point $u\in\R^{\mathcal{A}}$
   it defines a plane tropical curve $\mathcal{C}_F$ as the locus of
   non-differentiability of the tropical polynomial \label{page:troppoly}
   \begin{displaymath}
     F=\max\{u_{ij}+i\cdot x+j\cdot y\;|\;(i,j)\in\mathcal{A}\},
   \end{displaymath}
   and it defines a regular subdivision of $\Delta$.
   Each marked polygon of the subdivision is dual to a
   vertex of $\mathcal{C}_F$, and each edge $e$ of a marked polygon is dual
   to an edge $E$ of $\mathcal{C}_F$.
   Moreover, the edge $E$ is orthogonal to its dual edge $e$.
   Finally, the edge $E$ is unbounded if and
   only if its dual edge $e$ is contained in the boundary of the
   polygon $\Delta$.
   The \emph{weight} of an edge $E$ is equal to $\#(e\cap \Z^2)-1$.

   The duality implies that we can deduce the type of the marked
   subdivision from the plane tropical curve $\mathcal{C}_F$, but not the
   markings. To deduce the markings, we need to know the
   coefficients $u_{ij}$.

   Obviously, the vector $(1,\ldots,1)$ is contained in the lineality space of the secondary fan.
   Therefore we can mod out this vector and consider the resulting fan
   in $\R^{s-1}= \R^{\mathcal{A}}/(1,\ldots,1)\cdot \R$ with $s=\#\mathcal{A}$.
   We have seen above that every point in $\R^{\mathcal{A}}$
   defines a tropical curve via the tropical polynomial $
   \max\{u_{ij}+i\cdot x+j\cdot y\}$. Of course, adding $1$ to each
   coefficient $u_{ij}$ does not change the tropical curve associated
   to this tropical polynomial. Hence if we consider $\R^{\mathcal{A}}$ as
   a parametrising space for tropical curves, it makes sense to mod
   out $(1,\ldots,1)\cdot\R$, and we will do so in what follows. By abuse
   of notation, we call the fan in $\R^{s-1}$ that we get from the
   secondary fan in this way also the \emph{secondary fan}.

   \subsection{The dimension of cones in the secondary fan}\label{subsec-dimofsecfancones}

   Let $T=\{(Q_l,\mathcal{A}_l)\;|\;l=1,\ldots,k\}$ be a regular marked
   subdivision of $\Delta$. Let
   \begin{displaymath}
     L:=\Big\{(\lambda_{ij})\in \R^{\mathcal{A}} \;\Bigm|\; \sum\nolimits_{ij}
     \lambda_{ij}\cdot (i,j)=0, \;\sum\nolimits_{ij}\lambda_{ij}=0\Big\}
   \end{displaymath}
   be the space of affine relations among the lattice points $(i,j)$ of $\Delta$.
   For any $l$, let
   \begin{displaymath}
     L_{\mathcal{A}_l}= \{(\lambda_{ij})\in L\;|\;\lambda_{ij}=0
     \mbox{ for } (i,j)\notin \mathcal{A}_l\}
   \end{displaymath}
   be the space of affine relations among the elements of
   $\mathcal{A}_l$. Let $L_T$ be the sum $\sum_l L_{\mathcal{A}_l}$.

   \begin{lemma}
     The codimension of the cone of the secondary fan corresponding to
     the marked subdivision $T$ equals $\dim(L_T)$.

     In particular, a cone in the secondary fan corresponding to a
     marked subdivision is top-dimensional if and only if the marked
     subdivision is a triangulation, i.e. all polygons $Q_i$ are
     triangles and in each $Q_i$ no other point besides the vertices
     is marked.
   \end{lemma}
For a proof, see \cite[Corollary 2.7]{GKZ}.

   \begin{example}
     Let $T=\{(Q_1,\mathcal{A}_1),(Q_2,\mathcal{A}_2)\}$ be the
     subdivision shown in Figure \ref{fig:affinerelations},
     \begin{figure}[h]
       \centering
       \input{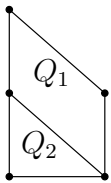}
       \caption{A marked subdivision}
       \label{fig:affinerelations}
     \end{figure}
     then $L$ is the kernel of the matrix
     \begin{displaymath}
       A=\left(
         \begin{array}[m]{ccccc}
           1&1&1&1&1\\0&1&0&1&1\\0&0&1&1&2
         \end{array}
       \right),
     \end{displaymath}
     where the second and third entry of each column of $A$ corresponds to
     the coordinates of a lattice point. Thus $L$ is generated by
     \begin{displaymath}
       (1,-1,-1,1,0) \;\;\;\mbox{ and }\;\;\;
       (0,1,0,-2,1).
     \end{displaymath}
     Then $L_{\mathcal{A}_2}$ is set of vectors in $L$ where the
     fourth and the fifth component vanish, and it is thus the zero
     space, while in $L_{\mathcal{A}_1}$ the first component has to
     vanish and it is thus generated by $(0,1,0,-2,1)$. We get
     \begin{displaymath}
       L_T=L_{\mathcal{A}_1}=(0,1,0,-2,1)\cdot\R
     \end{displaymath}
     and the codimension of the cone in the secondary fan
     corresponding to the marked subdivision $T$ is one.
   \end{example}

   \begin{remark}\label{rem-circuit}
     A cone in the secondary fan is of codimension one if and only if
     one of the $\mathcal{A}_i$ of  the corresponding
     marked subdivision contains a circuit, and this circuit is unique. Here, a \emph{circuit} is a set of
     lattice points that is affinely dependent but such that each
     proper subset is affinely independent.
     Figure \ref{fig:circuits} shows all types of circuits that can
     appear for point configurations in the plane together with some
     marked subdivisions of codimension one.
     \begin{figure}[h]
       \centering
       \input{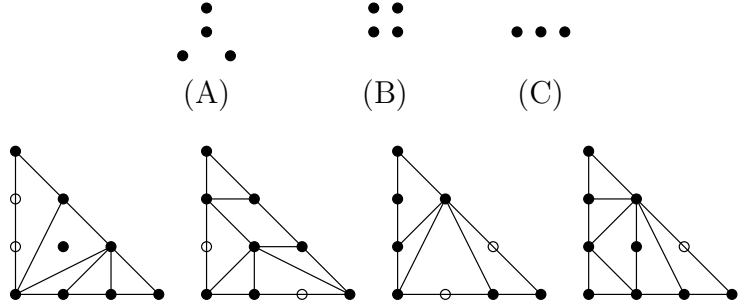}
       \caption{Planar circuits and marked subdivisions of codimension one}
       \label{fig:circuits}
     \end{figure}
   \end{remark}

   \subsection{The dimension of types of tropical curves}\label{subsec-dimofcurves}
   Given a tropical curve $C$, we have seen above that it is dual to a
   type $\alpha=\{Q_1,\ldots,Q_k\}$ of a marked subdivision. We call
   $\alpha$ also the \emph{type of the tropical curve}.
   We can parametrise all tropical curves of a given type by an
   unbounded polyhedron in $\R^{2+b}$ where $b$ denotes the number of
   bounded edges of $C$. This is true because we can move the curve in
   the plane, and we can change the lengths of the bounded edges
   without changing the type. However, the lengths cannot be changed
   independently if the tropical curve is of genus $g\geq 1$. We get $2g$
   (not necessarily independent) equations in $\R^{2+b}$ that tell us
   that the loops of $C$ have to close up. We define the
   \emph{type-dimension} $\tdim(\alpha)$ of a type $\alpha$ to be the
   dimension of the parametrising polyhedron.

   For the following lemma recall that we consider the secondary fan
   of $\Delta$ as a fan in $\R^{\mathcal{A}}/(1,\ldots,1)\cdot\R$.

   \begin{lemma}\label{lem-maxdim}
     Given a marked subdivision $T=\{(Q_l,\mathcal{A}_l)\}$ of
     $\Delta$ of type $\alpha$, we have
     \begin{displaymath}
       \tdim(\alpha) \leq  \dim(\sigma_T),
     \end{displaymath}
     where $\sigma_T$ denotes the cone of the secondary fan corresponding to $T$.

     Equality holds if and only if in $T$ all lattice points of
     $\Delta$ are marked, i.e.\ there are no white points.
   \end{lemma}
   \begin{proof}
     To a point $\overline{u}\in \sigma_T$ with representative
     $u\in\R^{\mathcal{A}}$ we associate a tropical polynomial
     $\max\{u_{ij}+ix+jy\}$ and thus a tropical curve of type
     $\alpha$. If we fix one of the polygons in $T$ and assign to the
     tropical curve the coordinates of the vertex corresponding to
     this polygon and the lattice lengths of the $b$ bounded edges, we
     get a map
     \begin{displaymath}
       \Phi_T:\sigma_T\longrightarrow \R^{2+b}
     \end{displaymath}
     from $\sigma_T$ to the parameter space of the type $\alpha$. This map
     is given by rational functions, since the coordinates of the vertices of the
     tropical curve and thus the lattice lengths of the edges are
     solutions of systems of linear equations of the form
     $u_{ij}+ix+jy=u_{kl}+kx+ly$.

     We have to show that $\Phi_S$ is a bijection onto the polyhedron
     parametrising the type $\alpha$ if $S$ is the subdivision of type
     $\alpha$ where all lattice points are marked.
     Then $\tdim(\alpha)=\dim(\sigma_S)$ in this case,
     and if $T$ is any other subdivision of type $\alpha$ then $\sigma_S$
     lies in the boundary of $\sigma_T$ and has strictly smaller
     dimension.

     Every tropical curve of type $\alpha$ comes from a point $u$
     (i.e. is the tropical curve associated to the tropical polynomial
     $\max\{u_{ij}+ix+jy\}$) which has to be inside a cone of the
     secondary fan corresponding to a marked subdivision of type
     $\alpha$.
     Assume now there is a lattice point which is not marked in this
     subdivision.
     That means that the corresponding term in the tropical polynomial
     can never be the maximum. Therefore we can vary the coefficient
     without changing the tropical curve, until it reaches the upper
     faces of the convex hull of $\{(i,j,u_{ij})\}\subset \R^3$. Thus
     every tropical curve of type $\alpha$ comes in fact from a point
     $u$ inside the cone $\sigma_S$ corresponding to the subdivision of type
     $\alpha$ where all lattice points are marked. This shows that
     $\Phi_S$ is surjective.

     In order to see that $\Phi_S$ is injective it suffices to show
     that the tropical curve defined by $\max\{u_{ij}+ix+jy\}$
     determines the class $\overline{u}$ of $u$ in $\sigma_S$ uniquely,
     since the polyhedron associated to the type
     $\alpha$ parametrises the tropical curves of type $\alpha$.

     The vertices of the tropical curve are the solutions of a system
     of linear equations of the form
     \begin{displaymath}
       \left(
         \begin{array}[m]{cc}
           i-k&j-l\\i-m&j-n
         \end{array}
       \right)
       \cdot
       \left(
         \begin{array}[m]{c}
           x\\y
         \end{array}
       \right)
       =
       \left(
         \begin{array}[m]{c}
           u_{kl}-u_{ij}\\u_{mn}-u_{ij}
         \end{array}
       \right)
     \end{displaymath}
     with an invertible coefficient matrix defined by an arbitrary
     choice of three vertices of the polygon dual to the vertex of interest
     in the tropical curve. The solution thus determines the
     inhomogeneity of the system, i.e.\ the $u_{ij}$ for $(i,j)$ a
     vertex in a fixed polygon $Q_k$ of $S$ are determined up to a
     common summand. Since each polygon $Q_i$ shares a vertex with
     some other $Q_j$, this shows that the $u_{ij}$ corresponding to
     vertices of some $Q_k$ are all determined up to a common summand.
     If $(i,j)$ is a lattice point in $Q_k$ which is not a vertex of
     $Q_k$, then $u_{ij}$ is determined by the $u_{mn}$ corresponding
     to vertices of $Q_k$, since $(i,j,u_{ij})$ is supposed to be
     visible in the upper face which projects to $Q_k$. But then
     $u_{ij}$ is determined up to the same summand as the
     $u_{m,n}$. Altogether this shows that $u$ is determined by the
     curve up to adding a multiple of $(1,\ldots,1)$, and the class
     $\overline{u}$ is determined uniquely.
   \end{proof}

   \section{The tropicalisation of the family of
     curves with a singularity in a fixed point}\label{sec-family}

   Fix a non-degenerate convex lattice polygon $\Delta$ and set
   $\mathcal{A}=\Delta\cap\Z^2=\{m_1,\ldots,m_s\}$.
   The closure of the image of the map
   \begin{displaymath}
     \Psi_{\mathcal{A}}: (\K^\ast)^2\rightarrow
     \PP^{s-1}_{\K}: (x,y)\mapsto (x^{m_{1,1}}y^{m_{1,2}},\ldots,x^{m_{s,1}}y^{m_{s,2}})
   \end{displaymath}
   is a  toric surface $\Sigma=\Tor_\K(\Delta)$ and the hyperplane
   sections are the closure of the images of the curves in
   $(\K^\ast)^2$ given by
   \begin{displaymath}
     f_{\underline{a}}=a_1 x^{m_{1,1}}y^{m_{1,2}}+\ldots+a_s x^{m_{s,1}}y^{m_{s,2}}=0,
   \end{displaymath}
   with $\underline{a}=(a_1,\ldots,a_s)$.

   The linear equations in the $a_i$ for the family
   $\Sing_\K(\Delta)$ of such curves with a singularity in the fixed point
   $\bp=(1,1)$ are
   \begin{displaymath}
     f_{\underline{a}}(1,1)=0,\;\;\;
     \frac{\partial f_{\underline{a}}}{\partial x} (1,1)=0,\;\;\;
     \frac{\partial f_{\underline{a}}}{\partial y} (1,1)=0,
   \end{displaymath}
   or equivalently we can say that the family $\Sing_\K(\Delta)$ is the
   kernel of the $3\times s$ matrix
   \begin{displaymath}
     A=\left(\begin{matrix}
         1&\ldots &1\\
         m_1&\ldots&m_s
       \end{matrix}\right).
   \end{displaymath}
   Notice that $A$ is just the matrix of our point configuration,
   after raising the points on the $\{t=1\}$-plane in $\R^3$, if we
   choose the coordinates $(t,x,y)$ on $\R^3$.

   We want to study the tropicalisation of $\ker(A)$,
   $$\Trop(\ker(A))=\Trop(\Sing_\K(\Delta)).$$

   \begin{remark}\label{rem-whichpoint}
     If we choose a different point $(p,q) \in (\C^\ast)^2$ and consider
     the family of curves with a singularity in
     $(p,q)$, then the coefficient matrix $A$ of the above linear equations changes.
     More precisely, it will be multiplied from the right by a
     diagonal matrix $D(p,q)=(d_{ij})_{i,j=1,\ldots,s}$
     with diagonal entry $d_{ii}=p^{m_{i1}}\cdot q^{m_{i2}}$ if $m_i=(m_{i1},m_{i2})$.
     Denote by $A(p,q)=A\cdot D(p,q)$ the new matrix, then the minors
     of $A(p,q)$ differ from the minors
     of $A$ only by certain factors, and each of these factors is a
     monomial in $p$ and $q$ since the
     columns in $A(p,q)$ differ from the corresponding columns of $A$
     only by factor which is a monomial
     in $p$ and $q$.

     However, the matroid of $A(p,q)$ is determined by the question which minors of $A(p,q)$
     vanish and which do not. So the matroids of $A$ and of $A(p,q)$ coincide.
     The tropical variety $\Trop(\ker(A))$ respectively $\Trop(\ker(A(p,q))$ depends only
     on the matroid of $A$ respectively of  $A(p,q)$
     (see \cite{Stu02}, \S\;9.3). Thus, the tropical variety of $\ker(A(p,q))$ is independent of
     the chosen point $(p,q) \in (\C^\ast)^2$.
   \end{remark}

We are grateful to Joaquim Ro\'{e} for pointing out the following:
\begin{remark}\label{rem-whichpoint2}
 For a point $(p,q)\in (\K^\ast)^2$, the family $\Sing_\K^{(p,q)}(\Delta)$ of curves with
   a singularity at the point $(p,q)$ tropicalises to a shift of the tropical variety $\Trop(\ker(A))$ we describe. 
This is true since a curve $f_{\underline{a}}$ is singular at $(p,q)$ if and only if $f_{ \Psi_{\mathcal{A}}(p,q)\cdot \underline{a}}$ is singular at $(1,1)$. Thus we have a map from $\Sing_\K^{(p,q)}(\Delta)$ to $\Sing_\K^{(1,1)}(\Delta)$ which sends $\underline{a}$ to $\Psi_{\mathcal{A}}(p,q)\cdot \underline{a}$, and the tropicalisation of this map is a shift.
\end{remark}

   \subsection{The tropicalisation of $\ker(A)$} \label{subsec-Bergman}

   Let us now study the tropicalisation of $\ker(A)$.  As remarked in
   the introduction, we consider this tropical variety as a fan
   in $\R^s/(1,\ldots,1)\cdot\R$.  By Section 2.5 of \cite{Spe05}, the
   fan is balanced.
   To study $\Trop(\ker(A))$, we use the following known results about
   the tropicalisation of linear spaces.

   It was observed in \cite{Stu02}, \S\;9.3, that the tropicalisation
   of a linear space $\ker(A)$ depends only on the matroid $M$
   associated to $\ker(A)$, as a linear subspace.
   This matroid can be specified by its collection of
   circuits, which are the minimal sets arising as supports of linear
   forms vanishing on $\ker(A)$ (resp.\ minimal sets arising as
   supports of elements in the row space of $A$).
   Equivalently, these are minimal sets $\{i_1,\ldots,i_r\}\subset
   \{1,\ldots,s\}$ such that the columns $b_{i_1},\ldots,b_{i_r}$ of a
   Gale dual $B$ of $A$ are linearly dependent.
   A \emph{Gale dual} is a matrix $B$ whose rows span the kernel of
   $A$.  Thus we can also describe the matroid associated to $\ker(A)$ as
   the matroid of the point configuration given by the columns of a
   Gale dual of $A$.

   Proposition 2.5 of \cite{FS05} states that the set of all $w\in\R^s$ such
   that $M_w$ (the matroid of bases $\sigma$ of $M$ for which $\sum_{i\in \sigma }w_i$ is maximal)
   contains no loop equals the tropicalisation of the
   linear space $\ker(A)$. In \cite{Stu02} and \cite{AK06} the set of all $w$ such that
   $M_w$ contains no loop is called the \emph{Bergman fan} of the
   matroid $M$. 
   Given $u\in \R^s$, let $\mathcal{F}(u)$ denote the unique \emph{flag of subsets}
   \begin{displaymath}
     \emptyset=: F_0 \subsetneqq F_1 \subsetneqq \ldots \subsetneqq
     F_k\subsetneqq F_{k+1}:= \{1,\ldots,s\}
   \end{displaymath}
   such that
   \begin{displaymath}
     u_i<u_j\;\;\;\Longleftrightarrow\;\;\;\exists\;m\;:\; i\in F_{m-1}
     \mbox{ and } j\not\in F_{m-1}.
   \end{displaymath}
   In particular,
   \begin{displaymath}
     u_i=u_j\;\;\;\Longleftrightarrow\;\;\;\exists\;m\;:\; i,j\in
     F_m\setminus F_{m-1}.
   \end{displaymath}
   The \emph{weight class} of a flag $\mathcal{F}$ is the set of all
   $u$ such that $\mathcal{F}(u)=\mathcal{F}$.
   We can describe weight classes by their defining equalities and
   inequalities.

   For example,  the set of all
   vectors $u$ satisfying $u_3<u_1=u_4<u_2$ defines a weight class in
   $\R^4$. It corresponds to the
   flag $\{3\} \subset \{1,3,4\} \subset \{1,2,3,4\}$.

   A flag $\mathcal{F}$ is a \emph{flag of flats} of the Gale dual $B$
   of $A$ respectively of the associated matroid $M$ if the linear span
   of the vectors $\{b_j \;|\; j\in F_i\}$ contains no $b_k$ with
   $k\notin F_i$. As before, the vectors $b_j$ denote the columns of a
   Gale dual of $A$.

   Theorem 1 of \cite{AK06} states that the Bergman fan of a matroid
   $M$ is the union of all weight classes of flags of flats of
   $M$. This result also follows from Theorem 4.1 of \cite{FS05}.

   As a consequence, we can study our tropical linear space by
   studying weight classes of flags of flats of a Gale dual of $A$.

   \begin{construction}\label{const-transf}
     Choose three points of $\mathcal{A}$ which are affinely independent.
     Then we can perform Gaussian elimination with the matrix $A$
     making the columns corresponding to these three points the columns
     with pivots.
     To the point configuration in threespace given by the columns of
     $A$, this Gaussian elimination has the effect of an affine
     transformation.
     Denote by $\widetilde{m}_i$ the $i$-th column of the transformed matrix $A$.
     To simplify notation, we will assume without restriction that the
     three points we chose are $m_1$, $m_2$ and $m_3$.
   \end{construction}

   \begin{remark}\label{rem-plane}
     Before we made the transformation from Construction \ref{const-transf},
     all columns of $A$ lived on the $\{t=1\}$ plane. Since  the three
     special points are transformed to
     $\widetilde{m}_1=(1,0,0)$, $\widetilde{m}_2=(0,1,0)$ and
     $\widetilde{m}_3=(0,0,1)$, the point configuration now sits on
     the $\{t+x+y=1\}$ plane spanned affinely
     by these three points.
   \end{remark}

   Performing this Gaussian elimination makes it easy to read off
   generators for the kernel of the matrix, and thus a possible Gale
   dual.

   \begin{example}\label{ex-galedual}
     Let  $\mathcal {A}$ be the point configuration in Figure \ref{fig:square}.
     \begin{figure}[h]
       \centering
       \input{./Graphics/square.pstex_t}
       \caption{A point configuration}
       \label{fig:square}
     \end{figure}
     Then $A$ is the matrix
     \begin{displaymath}A=\left(\begin{matrix}
           1&1&1&1&1&1&1&1\\
           0&1&0&1&2&0&1&2\\
           0&0&1&1&1&2&2&2
         \end{matrix}\right).\end{displaymath}
     We choose the first three points --- $(0,0)$, $(1,0)$ and $(0,1)$
     --- to be the three special points in Construction
     \ref{const-transf}. After performing Gaussian elimination, the
     matrix reads:
     \begin{displaymath}
       \widetilde{A}=\left(\begin{matrix}
           1&0&0&-1&-2&-1&-2&-3\\
           0&1&0&1&2&0&1&2\\
           0&0&1&1&1&2&2&2
         \end{matrix}\right)
       =\left(\mathbbm{1}_3\;|\;A_1\right).
     \end{displaymath}
     From this, we can easily read off a basis of the kernel.
     \begin{displaymath}
       B=\left(-A_1^t\;|\;\mathbbm{1}_{s-3}\right)=
       \left(\begin{matrix}
           1&-1&-1&1&0&0&0&0\\
           2&-2&-1&0&1&0&0&0\\
           1&0&-2&0&0&1&0&0\\
           2&-1&-2&0&0&0&1&0\\
           3&-2&-2&0&0&0&0&1
         \end{matrix}\right).
     \end{displaymath}
     Note that by construction, the negative of the first three
     entries of the $i$-th row are
     just the coordinates of the transformed point
     $\widetilde{m}_{i+3}$.
   \end{example}

   \begin{remark}\label{rem-galedual}
     Just as in Example \ref{ex-galedual} we have in general that the
     Gale dual $B$ constructed in this way has the following
     form:
     It is an $(s-3)\times s$-matrix where the first 3 (column)
     vectors $b_1$, $b_2$ and $b_3$ are the
     $(t,x,y)$-coordinates of the $(s-3)$ points $-\widetilde{m}_i$,
     $i=4,\ldots,s$, and the remaining vectors are the unit vectors
     $b_4=e_1,\ldots,b_s=e_{s-3}$.
     The column $b_i$ corresponds to the point $\widetilde{m}_i$ in the
     sense that the $i$-th entry of the rows of $B$ --- which are in
     the kernel of $\widetilde{A}$ --- gets multiplied with $\widetilde{m}_i$
     when computing the product $ \widetilde{A} \cdot B^t$.

     In such a Gale  dual we now want to find flags of flats, i.e.\
     flags of $s-3$ subspaces $V_i\subset \R^{s-3}$:
     \begin{displaymath}
       \{0\} \subsetneqq V_1 \subsetneqq \ldots \subsetneqq V_{s-3},
     \end{displaymath}
     where each $V_i$ is generated by a subset of the column vectors
     $b_j$ of the Gale dual indexed by the set $F_i$, and the vectors
     $\{b_j\;|\; j\in F_i\}$ are all the column vectors of the Gale dual
     that are contained in the subspace $V_i$.
     In particular, $F_{s-3}=\{1,\ldots,s\}$.
     We set $F_i':=F_i\setminus F_{i-1}$.
     Each $F_i'$ must of course consist of at least one element $j$.
     Since we have $s$ vectors in total, we have 3 ``extra'' vectors
     that can a priori belong to any of the $F_i'$.
     In the next lemma, we show that in fact we do not have that much choice.
   \end{remark}

   \begin{lemma}\label{lem-chains}
     With the notation of Remark \ref{rem-galedual}, for each flag of
     flats of a Gale dual $B$ of $A$ we have either
     \begin{enumerate}
     \item $\# F_i'=1$ for all $i=1,\ldots,s-4$ and $\# F_{s-3}'=4$, or
     \item $\# F_{s-3}'=3$ and there is a $j\in \{1,\ldots,s-4\}$ with $\# F_{j}=2$.
     \end{enumerate}
     In the first case, if $F_{s-3}'=\{a,b,c,d\}$, then any proper
     subset of the points $m_a$, $m_b$, $m_c$ and $m_d$ is affinely
     independent (i.e.\ $\{m_a,m_b,m_c,m_d\}$ is a circuit of type
     $(A)$ or $(B)$ as in Remark \ref{rem-circuit}).

     In the second case, if $F_{s-3}'=\{a,b,d\}$, the points $m_a$,
     $m_b$, $m_d$ are affinely dependent (i.e.\ $\{m_a,m_b,m_d\}$ is a
     circuit of type $(C)$ as in Remark
     \ref{rem-circuit}). Furthermore, all points $m_r$ with $r\in
     F_l'$, $l>j$, are on the same line as $m_a$, $m_b$ and $m_d$.
   \end{lemma}
   \begin{proof}
     First, we show that $\# F_{s-3}'$ cannot be 2 or 1. Assume it
     was.
     Then $\#F_{s-4}=s-2$ (resp.\ $s-1$) but the subspace $V_{s-4}$
     spanned by the vectors of $F_{s-4}$ is only $(s-4)$-dimensional.
     Remember that $b_4=e_1,\ldots,b_s=e_{s-3}$. To get an
     $(s-4)$-dimensional subspace with $s-2$ (resp.\ $s-1$) vectors,
     we in principal have 2 possibilities:
     \begin{enumerate}
     \item $\{b_r \;|\; r\in F_{s-4}\}$ can contain $s-4$ of the unit
       vectors $b_4,\ldots,b_s$, and 2 of the special vectors
       $b_1,b_2,b_3$, or
     \item it can contain $s-5$ (resp.\ $s-4$) unit vectors and all
       special vectors $b_1,b_2,b_3$.
     \end{enumerate}
     Let us consider case a) first. In $\{b_r \;|\; r\in F_{s-4}\}$,
     we are missing just one of the unit vectors, say $b_{j+3}=
     e_{j}$.
     Then the 2 special vectors which are also part of $\{b_r \;|\;
     r\in F_{s-4}\}$ must both have zeroes in the $j$-th component, or
     the dimension would be bigger than $s-4$. The $j$-th components
     of the two special vectors are two of the coordinates of the
     point $\widetilde{m}_{j+3}$. Thus without restriction
     $\widetilde{m}_{j+3}=(a,0,0)$ for some number $a$.
     But remember that the points $\widetilde{m}_i$ are all in the
     $\{t+x+y=1\}$ plane by Remark \ref{rem-plane}, thus $a=1$. But
     then $\widetilde{m}_{j+3}=\widetilde{m}_1=(1,0,0)$ which is a
     contradiction.

     Let us now consider case b). Now we are missing two of the unit
     vectors in $\{b_r \;|\; r\in F_{s-4}\}$, say $b_{j+3}=e_j$ and
     $b_{k+3}=e_k$ (resp.\ again just one). The three special vectors
     $b_1,b_2,b_3$ must have two linearly dependent rows in the $j$-th
     and $k$-th component, or the dimension would be bigger than $s-4$
     (resp.\ they must all have a zero in some row  which cannot be true
     since this row lives in the $\{t+x+y=1\}$ plane). But both rows again
     are the coordinates of the points $\widetilde{m}_{j+3}$ and
     $\widetilde{m}_{k+3}$. Both points live in the $\{t+x+y=1\}$ plane,
     so if they are linearly dependent, they are equal which is a
     contradiction.

     We conclude that for each flag of flats, $\# F_{s-3}'=4$ or $3$.
     Pick a flag of flats.
     For the statement about the affine dependencies, we want to
     switch to another Gale dual of $A$ which is more suitable for
     this particular flag of flats.
     First, choose two of the elements of $F_{s-3}'$, say $a$ and $b$,
     and an arbitrary $c$ such that $m_a$, $m_b$ and $m_c$ are
     affinely independent.
     We want to use these as the three pivot points in Construction
     \ref{const-transf}.
     We thus produce a new Gale dual (that must contain the same flag
     of flats). To make the notation simple, as before we want to call
     the special vectors $b_1$, $b_2$ and $b_3$, i.e. we assume
     without restriction that $a=1$, $b=2$ and $c=3$.
     Thus $1,2\in F_{s-3}'$.

     Let $d$ be a third index in $F_{s-3}'$ and assume that $m_1$, $m_2$ and $m_d$
     are affinely dependent.
     Note that affine dependence is preserved under the transformation
     that we perform to produce $\widetilde{A}$.
     We want to show that there cannot be a fourth element in
     $F_{s-3}'$, i.e.\ if $F_{s-3}'$ contains $a$, $b$ and $d$ with
     $m_a$, $m_b$ and $m_d$ affinely dependent, then $\# F_{s-3}'=3$.
     To see this, remember that the coordinates of the transformed
     point $\widetilde{m}_d$ appear in the $(d-3)$-rd row of the special
     vectors $b_1$, $b_2$ and $b_3$. But since $\widetilde{m}_d$,
     $\widetilde{m}_1=(1,0,0)$ and $\widetilde{m}_2=(0,1,0)$ are affinely
     dependent, it follows that $\widetilde{m}_d$ has a 0 as third
     coordinate. Thus $b_3$ has a $0$ in the $(d-3)$-rd row, and
     we have to show that
     $\{1,\ldots,s\}\setminus\{1,2,d,i\}\subseteq F_{s-4}$ implies
     $i\in F_{s-4}$.

     Assume first $i=3$. But $b_3$ is in the subspace generated by
     $\{b_4=e_1,\ldots,b_s=e_{s-3}\}\setminus \{b_d={e_{d-3}}\}$,
     since it has a $0$ in the $(d-3)$rd coordinate. Thus any set
     containing $\{1,\ldots,s\}\setminus\{1,2,d,3\}$ also contains
     $3$.

     Now assume $i\neq 3$. Suppose that $i$ is not in $F_{s-4}$. Then
     the $s-4$ vectors
     $\{b_3,b_4=e_1,\ldots,b_s=e_{s-3}\}\setminus \{b_i,
     b_d={e_{d-3}}\}$ generate the $s-4$-dimensional space
     $V_{s-4}$ and all vectors in $V_{s-4}$ have a zero in the
     $d-3$-rd component, since $b_3$ has so. Then the $i-3$-rd component
     of $b_3$ cannot be zero as well, since otherwise also the
     $i-3$-rd component of all vectors in $V_{s-4}$ would be zero in
     contradiction to the dimension being $s-4$. But then
     $b_i=e_{i-3}$ is a linear combination of
     $\{b_3,b_4=e_1,\ldots,b_s=e_{s-3}\}\setminus \{b_i,
     b_d={e_{d-3}}\}$ and it is in $V_{s-4}$, which implies that $i\in
     F_{s-4}$ in contradiction to our assumption. This shows that $i$
     is contained in $F_{s-4}$.

     Hence we have shown that
     $F_{s-4}=\{1,\ldots,s\}\setminus\{1,2,d\}$ and thus $\#
     F_{s-3}'=3$, if we assume that $1,2,d\in F_{s-3}'$ are affinely
     dependent.

     If $\# F_{s-3}'=4$, we thus have four points such that each
     proper subset is affinely independent, i.e.\ a circuit of type
     (A) or (B) from Remark \ref{rem-circuit}.

     Now assume that $F_{s-3}=\{a,b,d\}$. We have to show that $m_a$,
     $m_b$ and $m_d$ are affinely dependent, and that furthermore all
     points $m_r$ with $r\in F_l'$, $l>j$, are on the same line as
     $m_a$, $m_b$ and $m_d$, where $j$ is such that $\#F_j'=2$.
     Again, we want to pick a suitable Gale dual for this flag, i.e.\
     we assume without restriction that $a=1$ and $b=2$.
     As the third pivot point, we pick a point $m_c$ such that $m_a$,
     $m_b$ and $m_c$ are affinely independent, and such that $c$ is in
     $F_k'$ with $k$ maximal. That means, if $i\in F_l'$, $l>k$, then
     $m_a$, $m_b$ and $m_i$ are affinely dependent, i.e. $m_i$ lies on
     the line through $m_a$ and $m_b$. Again, we assume
     $c=3$. Now $F_{s-4}$ must contain all elements except $1$, $2$
     and $d$, thus $F_{s-4}=\{3,\ldots,s\}\setminus \{d\}$. $V_{s-4}$
     is an $s-4$-dimensional subspace. The vectors $\{b_4=e_1,\ldots
     b_s=e_{s-3}\}\setminus\{b_d=e_{d-3}\}$ thus span this subspace,
     and hence $b_3$ is a linear combination of these vectors. This
     implies that $b_3$ has a $0$ in the $(d-3)$-rd row, which in turn
     implies that $\widetilde{m}_d$, $\widetilde{m}_1=(1,0,0)$ and
     $\widetilde{m}_2=(0,1,0)$ are affinely dependent.

     Actually, $b_3$ has a $0$ in exactly those rows that correspond
     to points $\widetilde{m}_i$ which are affinely dependent of
     $\widetilde{m}_1=(1,0,0)$ and $\widetilde{m}_2=(0,1,0)$, i.e.\ that correspond
     to points $m_i$ on the same line as $m_1$ and $m_2$. If we set
     \begin{align*}
       B=&\{b_i\;|\; i\geq 4, m_i \mbox{ is not on the line through }m_1
       \mbox{ and }m_2\}\\
       =&\{e_{i-3}\;|\; i\geq 4, m_i \mbox{ is not on the line through }m_1
       \mbox{ and }m_2\},
     \end{align*}
     then $b_3$ is a linear combination of the elements of $B$ and non
     of the coefficients is zero. Thus any subset of $B\cup\{b_3\}$ of
     size $\#B$ is a basis of the span of $B\cup\{b_3\}$ which shows
     that some $V_l$ contains $\#B$ of the vectors of $B\cup\{b_3\}$ if
     and only if it contains all of them. Above we defined $k$ as the
     maximal index such that $F_k'$ contains an $i$ with $m_i$
     affinely independent of $m_1$ and $m_2$, then the previous
     considerations show that $V_k$ contains all vectors in
     $B\cup\{b_3\}$ while in $V_{k-1}$ two of them are missing. This
     shows that $F_k'$ has size two, i.e.\ $F_k'=F_j'$, (and contains
     the index $c=3$ by choice) and that for $l>k=j$ and $i\in F_l'$ the point $m_i$
     does not lie on the line through $m_1$ and $m_2$.
   \end{proof}

   \begin{remark}
     The following reversed statement of \ref{lem-chains} holds true as well:
     \begin{enumerate}
     \item For any circuit $\{m_a,m_b,m_c,m_d\}$ there exist all flags
       of flats satisfying $\#F_j'=1$ for all $j\neq s-4$ and
       $F_{s-4}'=\{a,b,c,d\}$.
     \item For any circuit $\{m_a,m_b,m_d\}$ and any choice of $m_c$
       and $m_e$ which are not on the line of $m_a$, $m_b$, $m_d$,
       there exist all flags of flats satisfying $F_{s-4}'=\{a,b,d\}$,
       $F_j'=\{c,e\}$, and all $i\in F_l'$ with $l>j$ satisfy $m_i$ is
       on the line.
     \end{enumerate}
     This can be seen similar to the proof of Lemma \ref{lem-chains}
     by picking a suitable Gale dual. In case (a), all vectors $b_i$
     with $i\notin F_{s-4}'$ are unit vectors and we can thus form any
     possible flag with them.
     In case (b), we can pick $a$, $b$ and $c$ as pivots, and then we
     can pick any flag such that $c$ and $e$ appear latest among all
     $i$ such that $m_i$ is not on the line of $m_a$ and $m_b$.
   \end{remark}

   \subsection{Steps towards the classification of tropical
     curves with a singularity in a fixed point}\label{subsec-class}

   As a consequence, we can try to classify all types of tropical
   curves with a singularity in a fixed point.
   To do this, let us first express the statement about the flags of
   flats from Lemma \ref{lem-chains} in terms of weight classes and
   marked subdivisions. We keep the notation from Remark
   \ref{rem-galedual}.
   The following list shows the important parts of the different
   weight classes we get and sums up what we can say about the marked
   subdivisions and their dual tropical curves.
   \begin{enumerate}
   \item \label{class-1}
     Assume we have a flag with $\#F_{s-4}'=4$ and the corresponding
     circuit is of type (A) or (B) as in Remark
     \ref{rem-circuit}. Then these points get the highest
     weight. Consequently, the triangle resp.\ quadrangle which is the
     convex hull of the circuit is part of the marked subdivision
     corresponding to any $u$ in the weight class.
     Besides, in the tropical polynomial which has $u$ as
     coefficients, the four terms corresponding to the four points
     have the same coefficients. The vertex dual to the triangle
     resp.\ quadrangle is at the point $(x,y)$ where the maximum is
     attained by those four terms, in particular the four terms are
     equal at this vertex. That means, we can set the four terms equal
     and solve for $x$ and $y$ to get the position of the
     vertex. But since the coefficients are all equal, we get $x=y=0$
     when solving.

     Thus the dual tropical curve has the point $\bx_0=(0,0)$ as a vertex of multiplicity
     strictly larger than one (corresponding to a triangle with an
     interior point, and thus of area bigger $1/2$), or it has a $4$-valent
     vertex at $(0,0)$ (see Figure~\ref{fig:type(a)}).
     \begin{figure}[h]
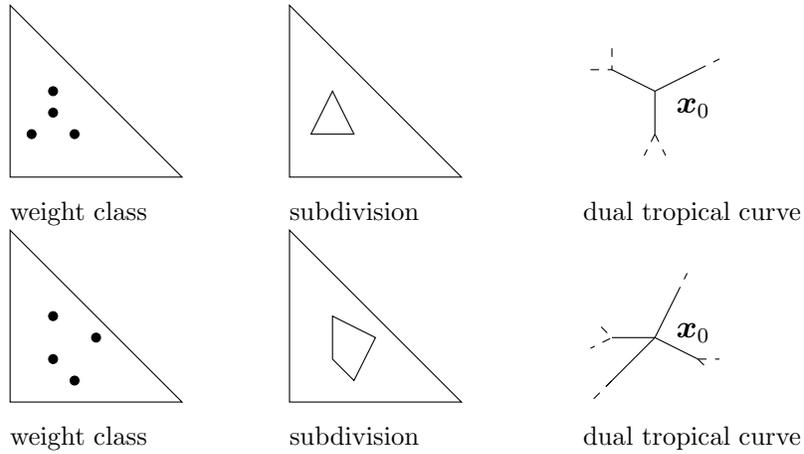

       \centering
       \input{./Graphics/weightclass.pstex_t}

       \input{./Graphics/weightclass2.pstex_t}
       \caption{Weight classes of type (a) and their tropical curves}
       \label{fig:type(a)}
     \end{figure}
   \item \label{class-3}Assume we have $F_{s-4}'=\{a,b,d\}$ and $F_j'=\{c,e\}$.
     In the picture, we draw the three points of highest weight black
     and the two points $m_c$ and $m_e$ in grey. Notice that the grey
     points are of the same height, namely the highest height of all
     points which are not on the line through $m_a$, $m_b$ and
     $m_d$. The points on this line can have higher heights however.

     Unfortunately, we cannot say much about the subdivision in this
     case. We can only be sure that the edge through
     $m_a$, $m_b$ and $m_d$ will be part of the subdivision. In the
     dual picture, this means we can see an edge of weight at least
     2. Furthermore, this edge must pass through the point
     $\bx_0=(0,0)$. The latter can be seen again by solving for the
     coordinates $(x_1,y_1)$ and $(x_2,y_2)$ of the two vertices
     adjacent to this edge. Since the heights of the three points
     $m_a$, $m_b$ and $m_d$ are equal, it follows that the line of
     which the dual edge is a segment passes through $(0,0)$. Since
     the height of any point which belongs to an adjacent polygon of
     the edge through $m_a$, $m_b$ and $m_d$ is below the height of
     these, it follows that $x_1<0$ and $x_2>0$ (or vice versa)
     (see Figure~\ref{fig:type(b)}).
     \begin{figure}[h]
       \centering
       \input{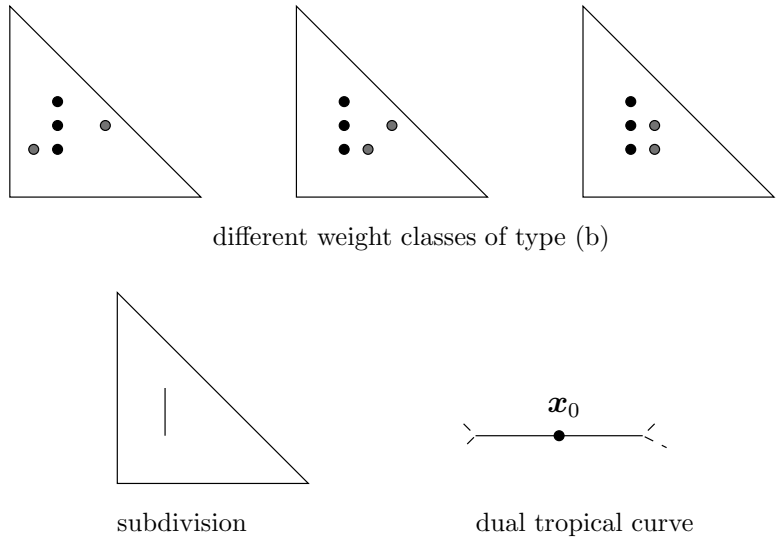}
       \caption{Different weight classes of type (b) and their tropical curve}
       \label{fig:type(b)}
     \end{figure}
   \end{enumerate}

   \begin{remark}
     The reason why we cannot say more than this is that we cannot
     predict how the polygons in the subdivision adjacent to the edge
     through $m_a$, $m_b$ and $m_d$ look like. It is possible that the
     grey points are not boundary points of an adjacent polygon. Even
     though they have the highest height of all points which are not
     on the line through $m_a$, $m_b$ and $m_d$, $(m_c,u_c)$ and
     $(m_e,u_e)$ could still lie
     below the upper faces of the convex hull of the points
     $(m_i,u_{i})$.
     As an example, take the point configuration in the picture below
     (where $m_c=(1,1)$) and take a weight vector $u$ as depicted in
     the middle.
     \begin{center}
       \input{./Graphics/exeat.pstex_t}
     \end{center}
     This weight vector is in the weight class
     \begin{displaymath}
       u_{h}<u_{f}<u_{g}<u_{c}=u_{e}<u_{a}=u_{b}=u_{d}
     \end{displaymath}
     which comes from the flag indexed by
     \begin{displaymath}
       \{h\} \subsetneqq  \{h,f\} \subsetneqq \{h,f,g\} \subsetneqq
       \{h,f,g,c,e\} \subsetneqq \{h,f,g,c,e,a,b,d\}.
     \end{displaymath}
     In the picture, we can see the marked subdivision induced by
     $u$. Note that the point $m_c$ is not part of a polygon adjacent
     to the edge through $m_a$, $m_b$ and $m_d$.
   \end{remark}

   A general point in a weightclass satisfies only the equalities given by
   the corresponding flag of flats, and strict inequalities otherwise. We can
   describe lower dimensional cones of the tropical variety $\Trop(\ker(A))$
   by forcing some of the inequalities to become equalities. Here, we
   restrict ourselves to the classification of top-dimensional cones.

   \subsection{$\Trop(\ker(A))$ and the secondary fan}\label{subsec-familyandsecfan}

   We have seen above that all the subdivisions we get in our family
   contain a circuit (either as a polygon $Q_i$, or as the face of a
   polygon $Q_i$). Hence the tropical variety $\Trop(\ker(A))$ lives
   inside the codimension-$1$-skeleton of the secondary
   fan.
   Furthermore, no weight class belonging to a flag of flats of type (a)
   of Subsection \ref{subsec-class} contains
   the lineality space of the secondary fan. The following lemma then shows that
   in a sense it is just the lineality space which is missing to
   pass from the cone of a weight class to a codimension one cone of
   the secondary fan.

   Remember that we have mod out
   the vector $(1,\ldots,1)$ already. But the secondary fan still
   contains a 2-dimensional \emph{lineality space} $L$ spanned by the
   vector consisting of the $x$-coordinates of the points $m_i$, and
   the vector consisting of their $y$-coordinates. This is true because
   if we incline the heights $u_i$ of the points $(m_i,u_{i})$ by a
   fixed multiple of the $x$-coordinates of the $m_i$ respectively of
   the $y$-coordinates of the $m_i$, we do not change the projection
   of the upper faces of the convex hull.

   \begin{lemma}\label{lem-squarecirc}
     Let $\Delta$ be a convex lattice polygon in the plane with associated
     matrix $A$ and Gale dual $B$ of $A$, and
     let $Z$ be a circuit in $\Delta$ of type (A) or (B) as in Remark
     \ref{rem-circuit}, i.e.\ a circuit consisting of four elements
     $Z=\{m_a,m_b,m_c,m_d\}$.

     Then the union of all weight classes
     $\tau_{\mathcal{F}}$ of flags of flats $\mathcal{F}$ of $B$ that end with
     $F_{s-4}'=\{a,b,c,d\}$ (where again we use the notation from
     \ref{rem-galedual}) plus the lineality space $L$ of the secondary
     fan of $\Delta$ equals the union of all codimension one cones $\sigma_T$ of the
     secondary fan of $\Delta$ corresponding to subdivisions $T$ that contain
     this circuit, i.e.
     \begin{displaymath}
       \left(\bigcup_{\mathcal{F}} \overline{\tau_{\mathcal{F}}}\right) +
       L =\bigcup_T \overline{\sigma_T},
     \end{displaymath}
     where the union on the left goes over all flags of flats
     $\mathcal{F}$ of $B$ that end with $F_{s-4}'=\{a,b,c,d\}$ and the union
     on the right goes over all subdivisions $T$ that contain the
     circuit $Z$.
   \end{lemma}
   \begin{proof}
     We have seen in our Classification \ref{subsec-class} that the
     marked subdivision of a vector $u$ in any weight class
     corresponding to such a flag of flats contains the circuit as a
     polygon. Thus $``\subset''$ is obvious.
     Pick any $u$ in $\sigma_T$, then we can write it as a sum of a vector
     in the lineality space and a vector that satisfies that the
     heights of the four points $m_a$, $m_b$, $m_c$ and $m_d$ are
     equal and highest among all heights. This shows $``\supset''$.
   \end{proof}

   Note also that the statement makes sense dimension-wise: The
   secondary fan is of dimension $s-1$ and the codimension 1 cone
   $\sigma_T$ of dimension $s-2$. Our tropical variety $\Trop(\ker(A))$ is
   of the same dimension as the ``classical'' variety $\ker(A)$ which
   is $s-4$-dimensional, since it lives in projective space of
   dimension $s-1$ and is given by $3$ independent equations.


   \begin{remark}\label{rem-special}
     Next we want to understand the cones of the secondary fan of
     $\Delta$ which correspond to flags of flats respectively weight
     classes of type (b) in the Classification \ref{subsec-class}. Let
     us thus assume that we have such a flag
     $\mathcal{F}$ of flats with $F_{s-3}'=\{a,b,d\}$ and $F_j'=\{c,e\}$ as in the
     proof of Lemma \ref{lem-chains}.

     The case where the
     points $m_c$ and $m_e$ span a line parallel to the line through
     $m_a$, $m_b$ and $m_d$, i.e.\ where we are in a situation as depicted
     in Figure \ref{fig:fuenferflaggen}, plays a special role.
     \begin{figure}[h]
       \centering
       \input{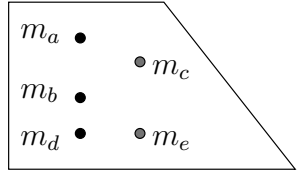}
       \caption{A weight class of type (b) in the boundary of others}
       \label{fig:fuenferflaggen}
     \end{figure}

     Let $\mathcal{F}=\mathcal{F}(u)$, let $T$ be the subdivision of
     $\Delta$ such that $u\in \sigma_T$ and let $Q$ be the polygon in $T$
     which contains the circuit $Z=\{m_a,m_b,m_d\}$ and lies on the
     same side of $Z$ as the points $m_c$ and $m_e$ (see Figure
     \ref{fig:5flaggen}). We then have to 
     distinguish two subcases.
     \begin{figure}[h]
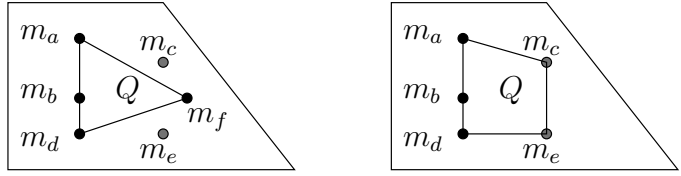

       \centering
       \input{./Graphics/fiveflag2.pstex_t}
       \hspace{1cm}
       \input{./Graphics/fiveflag1.pstex_t}
       \caption{Two different types of the boundary type}
       \label{fig:5flaggen}
     \end{figure}
     Either $Q$ contains a vertex whose
     distance to the line through $m_a$, $m_b$ and $m_d$ is larger
     than the distance of $m_c$ and $m_e$ to this line, or $Q$ is the
     polygon spanned by $Z$ and the two vertices $m_c$ and $m_e$ (see
     Figure \ref{fig:5flaggen}). To convince yourself that this is
     true recall that $u_a=u_b=u_d>u_c=u_e>u_i$ for all other $i$;
     thus, if $Q$ has no vertex whose distance is larger than that of
     $m_c$ and $m_e$ the planar polygon spanned by $(m_a,u_a)$, $(m_b,u_b)$,
     $(m_d,u_d)$, $(m_c,u_c)$ and $(m_e,u_e)$ in three space is an
     upper face of the extended Newton polytope corresponding to
     $\Delta$ and $u$.

     If $Q$ is spanned by $m_a,\ldots,m_e$ then the cone $\sigma_T$ is in
     the boundary of the cone $\sigma_S$ for a subdivision $S$ as shown in
     Figure  \ref{fig:5flagge}, where four of the lattice points form
     a quadrangle. Such quadrangles where already considered in Lemma
     \ref{lem-squarecirc}, and thus together with the cone $\sigma_T$ the
     weight class $\tau_{\mathcal{F}}\subset \sigma_T$ is contained in the
     boundary of cones of the secondary fan belonging to weight classes
     of type (a).
     \begin{figure}[h]
       \centering
       \input{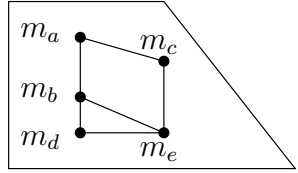}
       \caption{A subdivision such that $\sigma_S$ contains $\sigma_T$ in its boundary}
       \label{fig:5flagge}
     \end{figure}

     If $Q$ instead contains a vertex which is further away from the
     line through $m_a$, $m_b$ and $m_d$ than $m_c$ and $m_e$, then the
     cone $\sigma_T$ as well as the cone $\tau_{\mathcal{F}}\subset \sigma_T$ of
     the weight class $\mathcal{F}$ are in the boundary of cones of
     the secondary fan belonging to weight classes of type (b) as
     considered in the following Lemma \ref{lem-3ptscirc}.

     In any case it is not necessary to consider these weight classes
     in order to get a full picture of the codimension one cones of
     the secondary fan of $\Delta$ fixed by the weight classes of type
     (a) or (b).
   \end{remark}

   \begin{lemma}\label{lem-3ptscirc}
     Let $\Delta$ be a convex lattice polygon in the plane with associated
     matrix $A$ and Gale dual $B$ of $A$, and
     let $Z$ be a circuit of type (C) as in Remark \ref{rem-circuit},
     i.e.\ a circuit consisting of three elements
     $Z=\{m_a,m_b,m_d\}$.

     Then
     \begin{displaymath}
       \left(\bigcup_{\mathcal{F}} \overline{\tau_{\mathcal{F}} }\right)
       + L =\bigcup_T \overline{\sigma_T},
     \end{displaymath}
     where
     \begin{itemize}
     \item $L$ is the lineality space of the secondary fan of $\Delta$;
     \item the union on the left is the union of all weight classes
       $\tau_{\mathcal{F}}$ of flags of flats $\mathcal{F}$ of $B$ as in
       \ref{class-3} of the Classification \ref{subsec-class}, except
       for those considered in Remark \ref{rem-special}; that
       is, the flags end with $F_{s-4}'=\{a,b,d\}$, have
       $F_j'=\{c,e\}$ where the line through $m_c$ and $m_e$ is not
       parallel to the line through $m_a$, $m_b$ and $m_d$, while
       $m_i$ is on the latter line  for all $i\in F_l'$ for $l>j$;
     \item
       \begin{itemize}
       \item if $Z$ is not contained in the boundary of $\Delta$, the
         union on the right is the union of all codimension one
         cones $\sigma_T$ of the secondary fan of $\Delta$ that correspond to
         subdivisions $T$ containing $Z$;
       \item if $Z$ is contained in the boundary of $\Delta$, then the
         union on the right is the union of all codimension one
         cones $\sigma_T$ of the secondary fan of $\Delta$ that correspond to
         subdivisions $T$ containing $Z$, except those $T$ for which
         the triangle containing $Z$ has its third vertex at a point
         of minimal distance from $Z$.
       \end{itemize}
     \end{itemize}
   \end{lemma}

\begin{figure}[h]
     \centering
     \input{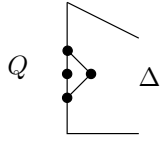}
     \caption{Triangulations we have to throw out}
     \label{fig:out}
   \end{figure}

   Figure \ref{fig:out} shows part of a triangulation corresponding  to one of
   the codimension one cones we throw out of the union if $Z$ is
   contained in the boundary of $\Delta$.

   \begin{proof}
     The proof is analogous to \ref{lem-squarecirc}.
     If $u$ is a vector in one of the weight classes we chose, then
     $u$ induces a subdivision containing $Z$. If $Z$ is in the
     boundary of $\Delta$, the two points $m_e$ and $m_c$ have to be
     on one side of the line through $Z$.
     We have to show that the polygon containing $Z$ is not a triangle
     with its third vertex at minimal distance.
     One of the points of $m_e$ and $m_c$ has to be at a non minimal
     distance, since we assume that they do not sit on a
     line parallel to the line through $Z$. Assume this point is
     $m_c$. One of the three lines
     connecting the points $(m_a,u_a)$ and $(m_c,u_c)$, resp.\
     $(m_b,u_b)$ and $(m_c,u_c)$, resp.\ $(m_d,u_d)$ and $(m_c,u_c)$,
     certainly lives above any line connecting $(m_a,u_a)$ with a
     point of minimal distance to $Z$. Hence $Z$ cannot be the face of
     a triangle with its third vertex at minimal distance.
     This proves ``$\subset$''.

     Conversely, we can write a vector $u$ in $\sigma_T$ as a sum of a
     vector in the lineality space and a vector that satisfies that
     the heights of the three points $m_a$, $m_b$ and $m_d$ are equal
     and highest, and that there are two points which are not on the
     line through $Z$ whose heights are equal and the
     highest among all points which are not on the line through $Z$.
     To do this, assume without restriction that $Z$ is on the line
     $\{x=0\}$. We can write $u$ as a sum of a multiple of the vector
     of $y$-coordinates of the $m_i$ and a vector $u'$ such that
     $u_a'=u_b'=u_d'$.

     \begin{figure}[h]
       \centering
       \input{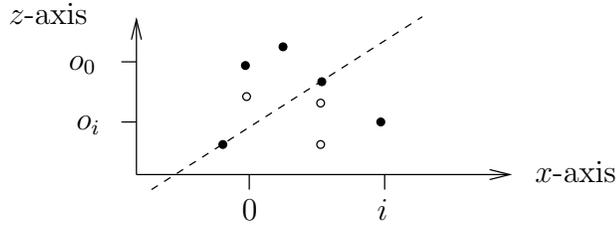}
       \caption{The projection of the $(m_i,u_i')$ to the $xz$-plane}
       \label{fig:theos}
     \end{figure}

     Now let $o_i$ be the maximum height on the line $\{x=i\}$, i.e.\
     $$o_i=\max\{u_j'\;|\; m_j\in \{x=i\}\},$$
     see Figure \ref{fig:theos}.
     Note that since $Z$ is part of the subdivision, $o_0$ must be the
     height of the points of $Z$.
     Even more, since $Z$ is part of the subdivision, the point
     $(0,o_0)$ will be above each line through two points $(k,u_i')$
     and $(l,u_j')$ with $m_i\in \{x=k\}$,
     $m_j\in \{x=l\}$ and $k<0<l$. This is true because otherwise
     there are two points in $Z$, say $m_a$ and $m_b$, which are not
     on different sides of the line through $m_i$ and $m_j$ and where
     one of them, say $m_a$, has a strictly larger distance to this
     line, and then the point $(m_b,u_b')$ would be strictly below the plane spanned by
     $(m_a,u_a')$, $(m_i,u_i')$ and $(m_j,u_j')$ in contradiction to
     the assumption that $Z$ is visible in the subdivision.
     Thus $o_0$ is contained in the boundary of the convex hull of the
     points $(k,o_k)$.
     Now we can add a multiple of the vector of $x$-coordinates to
     $u'$ to rotate the image in Figure \ref{fig:theos} about $(0,o_0)$ such
     that $o_0$ becomes the largest among the $o_k$ and such that
     the two next smaller $o_k$ have the same height, i.e such that
     there are $j\not=l$ with $o_j=o_l\geq o_h$
     for all $h\neq 0$.
     It is possible to make $o_0$ maximal, since $(0,o_0)$ is in the boundary
     of the convex hull of the points $(k,o_k)$.

     If $Z$ is not contained in the boundary of $\Delta$, we have
     points $(k,o_k)$ with positive and negative $k$-coordinate. By
     rotating about $(0,o_0)$ we can ensure that the vertices in the convex hull
     of the $(k,o_k)$ closest to the vertex $(0,o_0)$ on each side will have the
     same height, i.e.\ the two largest $o_k$ on each side of $Z$ will
     be equal. If
     $Z$ is contained in the boundary, we have only points $(k,o_k)$
     with positive $k$-coordinate (without restriction). However, the
     point of minimal distance $(1,o_1)$ is not a vertex of the convex
     hull of the points $(k,o_k)$. This is true since the triangle
     containing $Z$ does not have its vertex on the line
     $\{x=1\}$. This means again that we can make the two next largest
     heights $o_k$ equal by rotating. The point $u''$ we get in this way lives in a
     weight class as in \ref{class-3} of \ref{subsec-class}.
     This proves ``$\supset$''.

   \end{proof}

   \subsection{The connection to the tropical discriminant}\label{subsec-discriminant}
   The tropical discriminant has been studied by Dickenstein,
   Feichtner and Sturmfels (\cite{DFS05}).
   Their main result is that the tropicalisation of the discriminant
   of a point configuration $\mathcal{A}=\{m_1,\ldots,m_s\}$ --- i.e.\
   the locus of all parameters $\underline{a}$ for which the  curves
   $V(f_{\underline{a}})$ given by a polynomial
   \begin{displaymath}
     f_{\underline{a}}=a_1 x^{m_{1,1}}y^{m_{1,2}}+\ldots+a_s x^{m_{s,1}}y^{m_{s,2}}
   \end{displaymath}
   are singular --- equals $\Trop(\ker(A))+
   \mbox{rowspace}(A)$, where $+$ here denotes the Minkowski sum.
   This follows by a tropical version of Horn uniformisation. A curve
   $V(f_{\underline{a}})$ is singular in a point $(p,q)$ in the torus
   if and only if $V(f_{\Psi_{\mathcal{A}}(p,q)\cdot \underline{a}})$ (see
   beginning of Section \ref{sec-family}) is singular in $(1,1)$.
   This helps to express every point in the discriminant as the image
   under Horn uniformisation of a tuple consisting of a point in
   $\ker(A)$ and a point in the torus.
   Notice that the rowspace of $A$ equals the lineality space $L$ of
   the secondary fan of the point configuration.
   In the previous Section, we have described what cones of the
   secondary fan we get if we add this lineality space to our tropical
   variety $\Trop(\ker(A))$. Thus our result can also be seen as a
   description of the tropical discriminant.
   It follows that the tropical discriminant of a plane point
   configuration is a subfan of the secondary fan and consists of all
   closed codimension one cones of the secondary fan except the ones
   involving a circuit $Z$ consisting of three points on the boundary of
   $\Delta$ such that the triangle containing $Z$ has its third
   vertex at a point of minimal distance of $Z$.
   This description of the tropical discriminant was known before (see
   11.3.9 of \cite{GKZ}). There, $\Delta$-equivalent triangulations of
   the secondary fan are classified. Two triangulations are
   $\Delta$-equivalent, if their corresponding cones lie in the same
   top-dimensional cone of the Gr\"obner fan of the
   discriminant. Since the tropicalisation of the discriminant equals
   the codimension $1$-skeleton of the Gr\"obner fan, this means that
   two neighbouring triangulations are $\Delta$-equivalent if and only
   if the codimension $1$-cone that they meet in does not belong to
   the tropical discriminant.
   The only codimension $1$-cones of the secondary fan which do not
   belong to the tropical discriminant are the ones containing a
   circuit $Z$ of three points on the boundary such that the
   triangle containing $Z$ has its third vertex at a point of minimal
   distance of $Z$. Hence two neighbouring triangulations are
   $\Delta$-equivalent if and only if we can go from one to the other
   by a modification along such a circuit.

   \section{Classification of tropical curves of maximal
     dimensional type with a singularity in a fixed point} \label{sec-maxclass}

   We can say more if we restrict ourselves to tropical curves of
   maximal dimensional type.
   We have seen in Lemma \ref{lem-maxdim} that the dimension of a cone
   $\sigma_T$ of the secondary fan equals the dimension of its type if and
   only if the marked subdivision $T$ has no white points.
   Thus we can get tropical curves of maximal dimensional type only if
   we restrict ourselves to marked subdivisions corresponding to cones
   of smallest codimension and without any white points.
That is, we consider $\Trop(\Sing_\K(\Delta))$ with the stratification given by marked subdivisions, and we denote by $\MaxTrop(\Sing_\K(\Delta))$ the subcomplex consisting of cones corresponding to marked subdivisions for which all points are visible.
As remarked in the introduction, $\MaxTrop(\Sing_\K(\Delta))$ is the parameter space for tropical curves with a singularity in $(0,0)$.
By the following Lemma, $\MaxTrop(\Sing_\K(\Delta))$ is not a tropical
variety, i.e.\ it is not balanced.
\begin{lemma}
 The tropicalisation of a linear space is an irreducible tropical variety, i.e.\ contains no balanced proper subset of the same dimension.
\end{lemma}
\begin{proof}
Consider a codimension one cone $\tau$ of the tropicalisation of the linear space. By \cite{MS09}, Proposition 3.3.5, the star around $\tau$ is itself the tropicalisation of a linear ideal, in particular, it is a tropical line which is also a fan. Since tropical lines which are fans in $\R^n$ consist of $n+1$ rays of directions $-e_1,\ldots,-e_n$ and $e_1+\ldots,+e_n$, no proper subset of the top-dimensional neighbours of $\tau$ can satisfy a balancing condition. Thus the tropicalisation of a linear ideal is locally irreducible around each cone of codimension one.
By \cite{BJSST07}, Lemma 3.1, the tropicalisation of a linear ideal (and thus prime ideal) is connected in codimension one. It follows that the tropicalisation of a linear ideal is irreducible.
\end{proof}

\begin{figure}[h]
     \centering
     \input{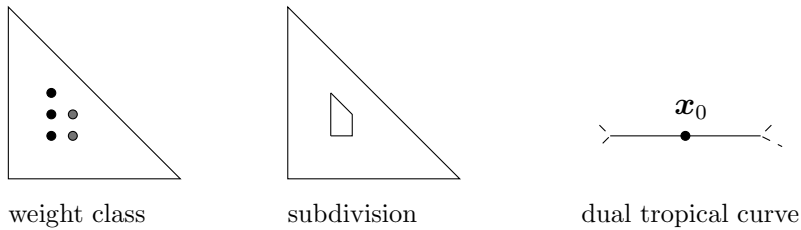}
     \caption{Subdivisions of codimension two}
     \label{fig:codim2}
   \end{figure}

 Note that in  Case \ref{class-3} of Classification
   \ref{subsec-class}, top-dimensional cones of $\Trop(\Sing_\K(\Delta))$ can
   also partly live in cones of codimension two of the secondary fan
   (see also Remark \ref{rem-special}).
   This is true because the two grey points $m_c$ and $m_e$ can be on
   a line which is parallel to the circuit $Z=\{m_a,m_b,m_d\}$. If
   these two points can be seen in the subdivision, then it belongs to
   a cone of the secondary fan of codimension two (see Figure
   \ref{fig:codim2}).

   In fact, we can relate these weight classes (where we restrict to
   the parts where the two points can be seen) to the secondary fan in
   a way similar to Lemma \ref{lem-squarecirc} and
   \ref{lem-3ptscirc}. We have to add only part of the lineality space
   of the secondary fan however. We have to add the vector consisting
   of all $y$-coordinates of the $m_i$ (if we assume without restriction that the
   circuit $Z$ is on the line $\{x=0\}$). On the right, we get the
   union over all codimension $2$ cones $\sigma_T$ of the secondary fan
   whose corresponding marked subdivision $T$ contains the polygon
   $\conv\{m_a,m_b,m_c,m_d,m_e\}$ and has all those points
   marked. This is true because for any vector $u\in \sigma_T$, we can add
   a multiple of the vector of $y$-coordinates of the $m_i$ to make the heights
   satisfy $u_a=u_b=u_d$ and $u_e=u_c$.

   In order to get tropical curves of maximal dimensional type, we
   thus have to study codimension $1$ cones of the secondary fan that
   are part of the tropical discriminant, and codimension $2$ cones
   that correspond to a marked subdivision containing a polygon
   $\conv\{m_a,m_b,m_c,m_d,m_e\}$ with all those points marked and
   such that $m_a$, $m_b$ and $m_d$ are on a line and $m_c$ and $m_e$
   are on a parallel line. We do not allow white points in the
   corresponding marked subdivisions.

   Now we classify top-dimensional cones of $\MaxTrop(\Sing_\K(\Delta))$.
   We go through the classification in \ref{subsec-class} and check
   what information on the dual tropical curve we can deduce in
   addition by assuming that there are no white points in the marked
   subdivision.
   \begin{enumerate}
   \item[(a)] Just as in \ref{subsec-class} \ref{class-1} we get tropical
     curves with a vertex of multiplicity $3$ at $\bx_0=(0,0)$ which is dual
     to a triangle with one interior lattice point, resp.\ with a
     $4$-valent vertex at $(0,0)$ whose dual polygon is a quadrangle
     not covering any other lattice points (see Figure \ref{fig:amax}).
     \begin{figure}[h]
       \centering
       \input{./Graphics/amax.pstex_t}
       \input{./Graphics/bmax.pstex_t}
       \caption{Types (a)}
       \label{fig:amax}
     \end{figure}

   \item[(b.1)] Let us consider a flag of flats as in \ref{subsec-class}
     \ref{class-3}. Since we do not want any white points, the two
     grey points have to be of minimal distance to the circuit $Z$,
     and they have to be vertices of polygons of the subdivision. The
     first case is that they are on different sides of $Z$ (see Figure
     \ref{fig:differentsides}).
     \begin{figure}[h]
       \centering
       \input{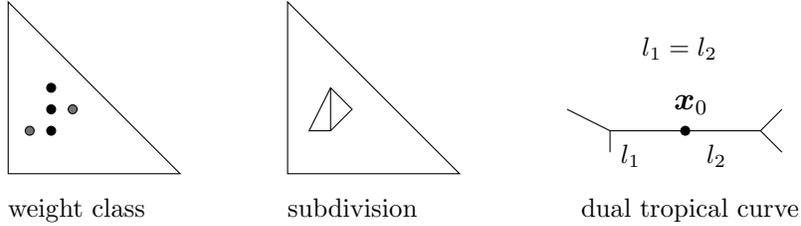}
       \caption{Type (b) with $m_c$ and $m_e$ on different sides}
       \label{fig:differentsides}
     \end{figure}
     Again, we can solve for the coordinates $(x_1,y_1)$ and
     $(x_2,y_2)$ of the two vertices adjacent to the edge through
     $(0,0)$. Now we know that dual to these vertices, we have two
     triangles whose third vertices are at the same height.
     If we assume without restriction that the circuit is on the line
     $\{x=1\}$ and the vertex of the left triangle is at $m_c=(0,0)$,
     then the equations to solve for $(x_1,y_1)$ and $(x_2,y_2)$ read:

     \begin{align*}
       &\lambda=\mu+x_1+(m_a)_2\cdot y_1=\mu+x_1+(m_b)_2\cdot y_1\\&
       \lambda+2\cdot x_2+(m_e)_2\cdot y_2= \mu+x_2+(m_a)_2\cdot
       y_2=\mu+x_2+(m_b)_2\cdot y_2
     \end{align*}
     where $\lambda$ is the height of the two grey points and $\mu$ is
     the height of the circuit points.
     Without restriction we can assume that $\lambda=0$ and $\mu>0$.
     Thus we conclude that the first vertex is at $(-\mu,0)$ and the
     second one is at $(\mu,0)$. In particular, the distances of both
     vertices to the singular point $(0,0)$ on the edge are
     equal.
   \item[(b.2)] Let us still consider flags as in
     \ref{subsec-class}\ref{class-3}, but now with the two grey points
     on the same side of the circuit $Z$. Again, the grey points have
     to be of minimal distance, and they have to be seen in the
     subdivision. Thus we can see a quadrangle with two parallel lines
     in the subdivision.
If $Z$ is not contained in the boundary of $\Delta$, there must be a
triangle whose vertex is of minimal distance in the subdivision on the
other side of $Z$ (see Figure
     \ref{fig:sameside}).
     \begin{figure}[h]
       \centering
       \input{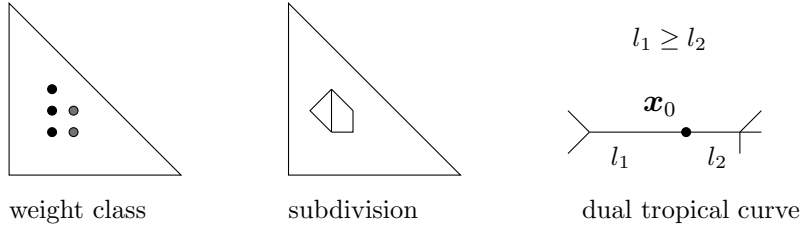}
       \caption{Type (b) with $m_c$ and $m_e$ on the same side}
       \label{fig:sameside}
     \end{figure}
     As above, we solve for the coordinates of the two vertices
     corresponding to the quadrangle and the triangle.
     Again, without restriction let us assume that the circuit is on
     the line $\{x=1\}$ and that the vertex of the left triangle is at
     $m_c=(0,0)$. Then the equations read:
     \begin{align*}
       &\nu=\mu+x_1+(m_a)_2\cdot y_1=\mu+x_1+(m_b)_2\cdot y_1\\&
       \lambda+2\cdot x_2+(m_e)_2\cdot y_2= \mu+x_2+(m_a)_2\cdot
       y_2=\mu+x_2+(m_b)_2\cdot y_2
     \end{align*}
     where $\lambda$ is the height of the grey points, $\mu$ is the
     height of the circuit points and $\nu$ is the height of the
     vertex of the left triangle. Without restriction, we can assume
     $\nu=0$, and $0<\lambda <\mu$.
     Thus the $3$-valent vertex is at $(-\mu,0)$ and the $4$-valent
     vertex is at $(\mu-\lambda,0)$.
     In particular, the distance from the $4$-valent vertex to the
     singular point $\bx_0=(0,0)$ is smaller than the distance of the
     $3$-valent vertex to $(0,0)$.
     If $Z$ is contained in the boundary, we see just the $4$-valent
     vertex, and $(0,0)$ lies on an infinite edge adjacent to this vertex.
   \end{enumerate}
%

We can sum up the results of this Section as follows:
\begin{theorem}\label{thm-class}
 Let $f_{a}=\sum_{(i,j)\in\mathcal{A}}a_{i,j}\cdot x^i\cdot
     y^j  \in \K[x,y]$ define a curve in the toric surface $\Tor_\K(\Delta)$ associated to $\Delta$ with a singularity in the point $(1,1)$.
Assume the dual marked subdivision defined by the coefficients $-\val(a_{i,j})$ of the tropical polynomial $\trop(f_a)$ is of maximal-dimensional type, i.e.\ all lattice points are marked (see Subsection \ref{subsec-dimofcurves}).
Then the tropicalisation locally around the singular point $(0,0)$ looks like one of the following cases.
The singular point $(0,0)$ is either
\begin{itemize}
 \item at a crossing of two edges of weight one,
\item or at a $3$-valent
   vertex of multiplicity $3$ adjacent to three edges of weight one (a nodal vertex),
\item or at the midpoint of an edge of weight $2$ connecting two $3$-valent vertices,
\item or in the interval from the $4$-valent vertex to the midpoint of an edge of weight $2$ connecting a $3$-valent and a $4$-valent vertex,
\item or on an infinite edge of weight $2$ whose end point is a
  $4$-valent vertex.
\end{itemize}

\end{theorem}

\begin{remark}
$\MaxTrop(\Sing_\K(\Delta))$ is connected in codimension one and to each cone of codimension one there are at least two top-dimensional cones adjacent.
To see this, consider two top-dimensional cones $\sigma_1$ and $\sigma_2$ of $\MaxTrop(\Sing_\K(\Delta))$. 
Both correspond to curves of maximal dimensional types $\alpha_1$ resp.\ $\alpha_2$, i.e.\ the dimension of the corresponding cones of their marked subdivisions in the secondary fan equals $\tdim(\alpha_i)$. We have seen that the corresponding cones in the secondary fan are either of codimension one or two, i.e.\ of dimension $s-2$ resp.\ $s-3$. Thus also $\tdim(\alpha_i)$ equals $s-2$ resp.\ $s-3$. Inside the polyhedra parametrising tropical curves of type $\alpha_i$ however, we only consider tropical curves with a singularity at $(0,0)$. For a type of dimension $s-2$, this means that we require a vertex resp.\ the midpoint of an edge to be at $(0,0)$. This yields two linear conditions and thus cuts out an $s-4$-dimensional part of the polyhedron. For a type of dimension $s-3$, we require $(0,0)$ to lie on an interval of a certain edge. This also cuts out an $s-4$-dimensional part.
It follows that $s-4$ is the maximal number of points in $\R^2$ (in general position) through which a tropical curve in $\MaxTrop(\Sing_\K(\Delta))$ can be fixed. That is, $s-4$-tuples of points in $\R^2$ through which a singular tropical curve of type $\alpha_i$ passes fill up a top-dimensional polyhedron in $\R^{2(s-4)}$, $s-4$-tuples of points through which a singular tropical curve of a lower-dimensional cone of $\MaxTrop(\Sing_\K(\Delta))$ passes fill up a lower-dimensional polyhedron accordingly. This yields a covering of $\R^{2(s-4)}$ by polyhedra of different dimensions. For each $\alpha_i$, we choose a configuration in $\R^{2(s-4)}$ living in the $s-4$-dimensional polyhedron consisting of $s-4$-tuples of points in $\R^2$ through which a singular tropical curve of type $\alpha_i$ passes, and we choose the two configurations such that they are connected by a line that avoids polyhedra of dimension smaller than $s-5$. This line can be lifted to a rational curve in $(\K^\ast)^{2(s-4)}$: we can cut out the line by $s-3$ hyperplanes each of which can be lifted to the variety of a binomial. The set of all binomials --- if chosen generically --- yields a toric and thus rational one-dimensional variety in $(\K^\ast)^{2(s-4)}$. We parametrise this curve thus producing an $s-4$-tuple of points in $(\K^\ast)^{2(s-4)}$ depending on a parameter. We require the curves in $\Sing_\K(\Delta)$ to pass through this $s-4$-tuple depending on a parameter, which yields a rational irreducible curve in $\Sing_\K(\Delta)$. This curve tropicalises to a connected tropical curve in $\MaxTrop(\Sing_\K(\Delta))$ which connects two curves in $\sigma_1$ and $\sigma_2$ and passes only through cones of codimension one.
\end{remark}

   \section{The tropicalisation of the family of curves with a
     singularity in a fixed point which is not a torus point} \label{sec-nottorus}

   In Remark \ref{rem-whichpoint} and \ref{rem-whichpoint2} we have seen that for any choice of
   singular point $(p,q)\in (\K^{\ast})^2$, we get the same
   tropicalisation for the family of curves with a singularity in
   $(p,q)$.
   What happens if we allow a point which is not in the torus, say
   $(p,q)=(1,0)$?
   The matrix $A$ given by the three equations $f_{\underline{a}}(1,0)=0$, $\frac{\partial
     f_{\underline{a}}(1,0)}{\partial x}(1,0)=0$ and $\frac{\partial
     f_{\underline{a}}(1,0)}{\partial y}(1,0)=0$ reads
   \begin{displaymath}
     A=\left(\begin{matrix}
         1&\ldots &1&0&\ldots&0&0&\ldots&0\\
         m_{1,1}&\ldots&m_{k,1}&0&\ldots&0&0&\ldots&0\\
         0&\ldots&0&1&\ldots&1&0&\ldots&0
       \end{matrix}\right),
   \end{displaymath}
   where we assume that the first block corresponds to the points
   $m_1,\ldots,m_k\in \mathcal{A}$ that satisfy $m_{i,2}=0$, the
   second block corresponds to the points with $m_{i,2}=1$, and the
   last block corresponds to the points with $m_{i,2}>1$.
   We assume that we have at least $3$ points with $m_{i,2}=0$ and at
   least $2$ points with $m_{i,2}=1$.

   Let us compute a Gale dual for this matrix as in Section \ref{subsec-Bergman}.
   \begin{example}\label{ex-point10}
     Assume the point configuration is as depicted in Figure \ref{fig:point10}.
     \begin{figure}[h]
       \centering
       \input{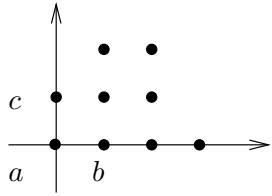}
       \caption{A point configuration}
       \label{fig:point10}
     \end{figure}
     Then the matrix A reads
     \begin{displaymath}
       A=\left(\begin{matrix}
           1&1&1&1&0&0&0&0&0\\
           0&1&2&3&0&0&0&0&0\\
           0&0&0&0&1&1&1&0&0
         \end{matrix}\right).
     \end{displaymath}
     Choose $a$, $b$ and $c$ as pivots and switch two columns such
     that the column of $c$ becomes the third column. Then the reduced
     row echelon form reads
     \begin{displaymath}
       \left(\begin{matrix}
           1&0&0&-1&-2&0&0&0&0\\
           0&1&0&2&3&0&0&0&0\\
           0&0&1&0&0&1&1&1&0
         \end{matrix}\right).
     \end{displaymath}
     The Gale dual we can easily read off from this form is
     \begin{displaymath}
       B=\left(\begin{matrix}
           1&-2&0&1&0&0&0&0&0\\
           2&-3&0&0&1&0&0&0&0\\
           0&0&-1&0&0&1&0&0&0\\
           0&0&-1&0&0&0&1&0&0\\
           0&0&-1&0&0&0&0&1&0\\
           0&0&0&0&0&0&0&0&1
         \end{matrix}\right).
     \end{displaymath}
     For the flags of flats, we can deduce the following:
     \begin{itemize}
     \item The vectors corresponding to points which are not on the
       line $\{y=0\}$ or $\{y=1\}$ (i.e.\ corresponding to the last
       block) are independent. They can be anywhere in a flag.
     \item For the vectors corresponding to the second block, i.e.\ to
       points on $\{y=1\}$: if (and only if) in the flag we collected
       all but one of those vectors, then the last one belongs to the
       subspace, too.
     \item For the vectors corresponding to the first block, i.e.\ to
       points on $\{y=0\}$: if (and only if) in the flag we collected all
       but two of those vectors, the last two belong to the subspace, too.
     \end{itemize}
     For the corresponding weight classes, we conclude:
     \begin{itemize}
     \item The maximum of heights appearing on the line $\{y=1\}$ is attained twice, and
     \item the maximum of heights appearing on the line $\{y=0\}$ is attained three times.
     \end{itemize}
     However, those two maxima can be in any relation to each
     other, and also in any relation to the heights of the points with
     $y$-coordinate larger than one.

     What can we conclude for the possible subdivisions?
     The only thing we know for sure is that the three points on
     $\{y=0\}$ must be seen in the subdivision.
     If they form a polygon with vertices in $\{y=1\}$, then these
     vertices have to be the two maximal points on this line.
     But they do not have to form a polygon with vertices in
     $\{y=1\}$, they could also form a polygon with a vertex in the
     line $\{y=2\}$.
     In Figure \ref{fig:point102}, we show several possible
     subdivisions. The three maximal points on $\{y=0\}$ are drawn in
     dark grey, the two maximal points on $\{y=1\}$ in light grey ---
     not depending on whether they can be seen in the subdivision or
     not.

     \begin{figure}[h]
       \centering
       \input{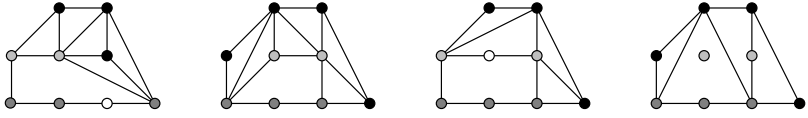}
       \caption{Some possible subdivisions}
       \label{fig:point102}
     \end{figure}

     For the dual tropical curves, we can conclude that there is an
     end of weight at least $2$ contained in the line $\{x=0\}$. (The
     $x$-coordinate can be found by solving the system of linear
     equations given by the three points in the subdivision. The fact
     that they are of the same height implies that the $x$-coordinate
     is $0$.)
     This fat end is either adjacent to an at least $4$-valent vertex,
     or to a $3$-valent vertex of multiplicity at least $4$.
     If the tropical curve is of maximal dimension, it has to end at a
     $4$-valent vertex.
     Since the negative of the valuation of the singular point $(1,0)$
     is $(0,-\infty)$, we expect the ``singularity information'' of
     the tropical curve to be contained in the ends.
     Local pictures of dual tropical curves are shown in Figure
     \ref{fig:point103}.

     \begin{figure}[h]
       \centering
       \input{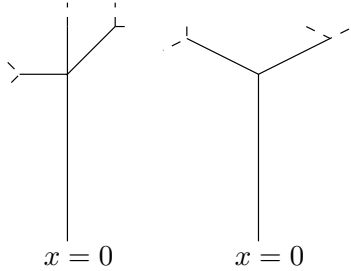}
       \caption{Dual tropical curves}
       \label{fig:point103}
     \end{figure}
   \end{example}

   \begin{remark}
     It is easy to show that the essential features of example
     \ref{ex-point10} hold in general. The matrix always consists of
     $3$ blocks corresponding to points in $\{y=0\}$, points in
     $\{y=1\}$ and points with $y>1$. The vectors in the Gale dual
     corresponding to the three blocks always behave like in the
     example, and we always get weight classes for the flags of flats
     where the maximal height in $\{y=0\}$ is attained three times and
     the maximal height in $\{y=1\}$ is attained twice. 
 \end{remark}
Thus, we can deduce the following result:
\begin{proposition}
 Let $f_{a}=\sum_{(i,j)\in\mathcal{A}}a_{i,j}\cdot x^i\cdot
     y^j  \in \K[x,y]$ define a curve in the toric surface $\Tor_\K(\Delta)$ associated to $\Delta$ with a singularity in the point $(1,0)$.
Then the tropicalisation has a
     fat end at $\{x=0\}$ (i.e.\ starting at the singular point $(0,-\infty)$) with either an at least $4$-valent vertex or
     a $3$-valent vertex of multiplicity at least $4$.
Assume the dual marked subdivision defined by the coefficients $-\val(a_{i,j})$ of the tropical polynomial $\trop(f_a)$ is of maximal-dimensional type, i.e.\ all lattice points are marked (see Subsection \ref{subsec-dimofcurves}).
Then the tropicalisation has a $4$-valent vertex
     adjacent to the fat end at $(0,-\infty)$.

\end{proposition}

     If we compare this to tropical curves with a singularity in
     $\bx_0=(0,0)$ that we studied in Section \ref{subsec-class}, we can see
     that we only get subdivisions where the circuit is of type $(C)$
     as in Remark
     \ref{rem-circuit} and is contained in the boundary on the line $\{y=0\}$.

   \section{Algebraic lifts of tropical curves of maximal
     dimensional type with a singularity in a fixed point} \label{sec-alglift}

   For the following considerations we assume that
   $\K=\overline{\C(t)}$ is the algebraic closure of the field of
   rational functions over the complex numbers.

   We would like
   to describe algebraic curves $C\in\Sing_\K(\Delta)$ which correspond
   to  tropical curves $\mathcal{C}\in\Trop(\Sing_\K(\Delta))$ of maximal
   dimensional type. We furthermore assume the following
   generality condition:
   \begin{enumerate}
   \item[(G)] $\mathcal{C}$ is a generic
     member in the interior of a top-dimensional cone of 
     $\MaxTrop(\Sing_\K(\Delta))$  and $C$ is a generic element of $\Sing_\K(\Delta)$
     with $\Trop(C)=\mathcal{C}$.
   \end{enumerate}
   Below we  specify this generality assumption
   which breaks certain explicit relations.

   As a particular consequence of our consideration, we give a
   conceptual explanation of the metric conditions for the type (b)
   curves (cf. Section \ref{sec-maxclass}).

 \subsection{ Tropical limits of plane algebraic curves over $\K$.} We
   shortly recall the definition of tropical limits
   used in the sequel following \cite{Shu06a}. A tropical curve $\mathcal{C}$ uniquely determines a
   convex piece-wise linear function $\nu:\Delta\to\R$ with
   $\max\nu=0$. Note that with the notation introduced in Section
   \ref{sec-secfan} on Page \pageref{page:troppoly} $\nu$ determines
   a defining tropical polynomial
   $F=\max\{u_{ij}+ix+jy\;|\;(i,j)\in\Delta\cap\Z^2\}$ for $\mathcal{C}$ via $\nu(i,j)=u_{ij}$.
   Without loss of generality, assume that $\mathcal{C}$ is defined over $\Q$ and
   that $\nu(\Delta\cap\Z^2)\subset\Z$ (the latter can be achieved by a
   suitable stretching of $\mathcal{C}$). An algebraic curve $C\in|{\mathcal
     L}_\Delta|$ with $\Trop(C)=\mathcal{C}$ is then given by an equation
   \begin{equation}
     f(x,y)=\sum_{(i,j)\in\Delta\cap\Z^2}(a^0_{ij}+O(t))t^{-\nu(i,j)}x^iy^j=0,\quad
     a_{ij}(t)\in\K\ ,\label{e2}
   \end{equation}
   where the $a_{ij}^0\in\C$ do
   not vanish since $(i,j)$ is visible in the subdivision $S_\mathcal{C}$ of $\Delta$
   induced by $\nu$, and the $O(t)$ are analytic functions in the disc
   $D_\varepsilon=\{|t|<\varepsilon\}$. Evaluating (\ref{e2}) for $t\in
   D_\varepsilon\backslash\{0\}$, we obtain a family of curves
   $C^{(t)}\subset\Tor_\C(\Delta)$ which admits a flat extension to $t=0$
   in the form
   $$\begin{matrix} {} &{}&{}&{}&\Tor(\widetilde\Delta)&{}&{}\\
     {} &{}&{}&{}&\parallel&{}&{}\\ C^{(0)}&\hookrightarrow&C & \hookrightarrow &
     \widetilde\Sigma&\hookleftarrow&\Sigma^{(0)}\\
     \downarrow&{}&\downarrow & {} & \downarrow & {} & \downarrow\\
     0&\in&D_\varepsilon & = & D_\varepsilon&\ni&0\end{matrix}$$ where
   $\widetilde\Delta=\{(\omega,z)\in\R^3\ \big|\ \omega\in\Delta,\
   z\le\nu(\omega)\}$ is the undergraph of $\nu$,
   $\Sigma^{(0)}=\bigcup_{\delta}\Tor(\delta)$ with $\delta$ ranging over all
   polygons of the subdivision $S_\mathcal{C}$, and $C^{(0)}\subset\Sigma^{(0)}$
   splits into the {\it limit curves} $C_\delta\subset\Tor(\delta)$ given
   by the equations
   $$f_\delta(x,y)\equiv\sum_{(i,j)\in\delta\cap\Z^2}a^0_{ij}x^iy^j=0\ .$$
   The data $(\mathcal{C},\{C_\delta\}_{\delta\in S_\mathcal{C}})$ is called the {\it tropical
     limit} of $C$.

   The fact that $C$ has a singularity at the point $\bp=(1,1)$ is
   equivalent to the fact that $C^{(t)}$ has a singularity at $(1,1)$
   for each $t\in D_\varepsilon\backslash\{0\}$ (cf. \cite[Lemma 2.3]{Shu06a}).
   The (constant) family of points
   $(1,1)\in(\C^*)^2\subset\Tor(\Delta)=\Sigma^{(t)}$, $t\in
   D_\varepsilon\backslash\{0\}$, has a limit point $p\in\Sigma^{(0)}$. Two
   cases are possible:
   \begin{enumerate}
   \item[(i)]
     $p\in(\C^*)^2\subset\Tor(\delta)$ for some polygon $\delta$ in
     $S_\mathcal{C}$; then, in particular, the equality $f(p)=f(1,1)=0$ implies
     that the initial form $\ini_{(0,0)}f(1,1)$ vanishes which forces
     the constancy of $\nu$ along $\delta$, and
     hence $\nu$ vanishes there by our assumptions; furthermore, the limit
     curve $C_\delta$ has a singularity at $(1,1)$; \item[(ii)]
     $p\in\Tor(\sigma)$, where $\sigma=\delta'\cap\delta''$ is a common side of
     polygons $\delta',\delta''$ in $S_\mathcal{C}$,
     $\Tor(\sigma)=\Tor(\delta')\cap\Tor(\delta'')$ is a common toric divisor;
     then $\nu$ vanishes along $\sigma$, and $p\in C_{\delta'}\cap
     C_{\delta''}\cap\Tor(\sigma)$, where the pairwise intersection
     multiplicities are $\ge 2$ (a transverse intersection point with
     $\Tor(\sigma)$ smoothes out in the deformation $C^{(0)}\to C^{(t)}$,
     $t\ne 0$, cf. \cite[Lemma 3.2]{Shu06a}).
   \end{enumerate} In the second case we shall {\it refine} the tropical
   limit of $C$ as described in \cite[Section 3.5 and 3.6]{Shu06a}.

   \medskip

   \subsection{ Curves of type (a).} 
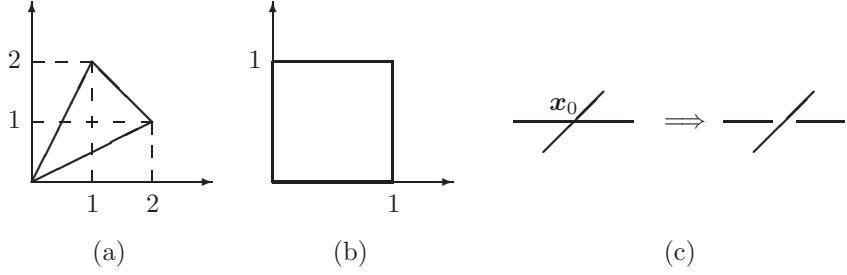
\begin{figure}[h]
     \setlength{\unitlength}{0.8cm}
     \begin{picture}(14,4.5)(0,0)
       \thinlines\put(0.5,1.5){\vector(0,1){3}}\put(0.5,1.5){\vector(1,0){3}}
       \put(4.5,1.5){\vector(0,1){3}}\put(4.5,1.5){\vector(1,0){3}}
       \dashline{0.2}(0.5,2.5)(2.5,2.5)\dashline{0.2}(0.5,3.5)(1.5,3.5)
       \dashline{0.2}(1.5,1.5)(1.5,3.5)\dashline{0.2}(2.5,1.5)(2.5,2.5)
       \thicklines\put(0.5,1.5){\line(1,2){1}}\put(0.5,1.5){\line(2,1){2}}
       \put(2.5,2.5){\line(-1,1){1}}\put(4.5,1.5){\line(1,0){2}}
       \put(4.5,1.5){\line(0,1){2}}\put(4.5,3.5){\line(1,0){2}}
       \put(6.5,1.5){\line(0,1){2}}\put(8.5,2.5){\line(1,0){2}}
       \put(9,2){\line(1,1){1}}\put(12.5,2){\line(1,1){1}}
       \put(12,2.5){\line(1,0){0.8}}\put(13.2,2.5){\line(1,0){0.8}}
       \put(0.1,2.4){$1$}\put(0.1,3.4){$2$}\put(4.1,3.4){$1$}
       \put(1.4,1){$1$}\put(2.4,1){$2$}\put(6.4,1){$1$}\put(11,2.4){$\Longrightarrow$}
       \put(1.5,0.2){\rm (a)}\put(5.5,0.2){\rm (b)}\put(11,0.2){\rm (c)}
       \put(9.1,2.7){$\bx_0$}
     \end{picture}
     \caption{Curves of type (a)}\label{fn1}
   \end{figure}

Let $\mathcal{C}$ be of type (a) introduced in
   Section \ref{sec-maxclass}. Then \begin{itemize}\item either the dual
     subdivision $S_\mathcal{C}$ of $\Delta$ consists of a triangle $\delta_0$ of
     lattice area $3$, which up to $SL(2,\Z)$-action and translations
     coincides with the triangle $\conv\{(0,0),(2,1),(1,2)\}$ (cf. Figure
     \ref{fn1}(a)), and the remaining pieces are primitive lattice
     triangles (i.e. of unit lattice area); \item or $S_\mathcal{C}$ contains a
     quadrangle $\delta_0$, which up to $SL(2,\Z)$-action and translations
     coincides with the square $\conv\{(0,0),(1,0),(0,1),(1,1)\}$ (cf.
     Figure \ref{fn1}(b)), and the remaining pieces again are primitive
     lattice triangles.\end{itemize} Observe that in this case all the
   edges $\sigma$ of the subdivision $S_\mathcal{C}$ have unit lattice length, and
   hence the limit curves $C_\delta$ can intersect the toric divisors
   $\Tor(\sigma)\subset\Tor(\delta)$ only transversally, which allows only
   the option (i) for the limit singular point $p\in\Sigma^{(0)}$
   described above. More precisely, the curve $C_{\delta_0}$ has a node
   at $\bp=(1,1)$ being irreducible if $\delta_0$ is a triangle (since the sides
   of $\Delta_0$ have unit length), or
   reducible if $\delta_0$ is a parallelogram, whereas the remaining
   limit curves (corresponding to primitive triangles) are
   nonsingular.
   In both the cases, $C$ is a nodal curve of genus
   $g=\#(\Int(\Delta)\cap\Z^2)-1$, and we can define a natural
   parametrisation of $\mathcal{C}$ of genus $g$ (cf. with a canonical
   tropicalisation in \cite{Tyo09}):
   \begin{itemize}\item If $\delta_0$ is a
     triangle, then $\mathcal{C}$ is self-parameterising of genus $g$. The singular point is a nodal vertex of the parametrising graph, since the image is dual to a polygon with one interior point. \item If
     $\delta_0$ is a parallelogram, then we resolve the $4$-valent vertex
     $\bx_0=(0,0)$ of $\mathcal{C}$, obtaining the refined tropical curve $\widehat {\mathcal{C}}$
     of genus $g$ as a parameterising graph of $\mathcal{C}$ (cf. Figure \ref{fn1}(c)). The singular point lifts to two interior points of edges of the parametrising graph.
   \end{itemize}

   \medskip

   \subsection{ Curves of type (b.1).} If $\mathcal{C}$ is of type (b.1) (see Section
   \ref{sec-maxclass} and Figure \ref{fig:differentsides}), then $S_\mathcal{C}$ consists of triangles,
   all of them but two $\delta',\delta''$ shown in Figure \ref{fig:differentsides} being
   primitive, and all the limit curves $C_\delta$ are nonsingular. Hence,
   we have the option (ii) for the limit singular point $p$, which then
   must belong to $\Tor(\sigma)$, where $\sigma=\delta'\cap\delta''$ is the
   common edge of length $2$ spanned by the circuit $Z$ (cf. Section
   \ref{sec-maxclass}). Furthermore, the limit curves $C_{\delta'}$,
   $C_{\delta''}$ must be quadratically tangent to $\Tor(\sigma)$ at
   $p$ (formally, we have an option of two
   transverse intersection points, however, it is not possible in
   our situation, since such points are smoothed out in the
   deformation $C^{(t)}$, $t\in(\C,0)$). 

\begin{figure}[h]
     \setlength{\unitlength}{0.83cm}
     \begin{picture}(14,8)(0,0)
       \thinlines\put(1,4.5){\vector(0,1){3.5}}\put(1,4.5){\vector(1,0){4}}
       \put(8,4.5){\vector(0,1){3}}\put(8,4.5){\vector(1,0){4}}
       \dashline{0.2}(8,6.5)(9,6.5)\dashline{0.2}(10,4.5)(10,6.5)
       \dashline{0.2}(0,2)(6,2)\dashline{0.2}(7.5,2)(12.5,2)
       \thicklines\put(0,1.5){\line(2,1){1}}\put(1,2){\line(0,1){1}}
       \put(1,2){\line(1,0){3}}\put(4,2){\line(1,1){1}}
       \put(4,2){\line(1,-1){1}}\put(8.5,2){\line(2,-1){1.5}}
       \put(11.5,2){\line(-2,-1){1.5}}\put(10,0){\line(0,1){1.25}}
       \put(3,5.5){\line(-1,2){1}}\put(2,7.5){\line(1,0){1}}
       \put(3,5.5){\line(0,1){2}}\put(3,5.5){\line(1,1){1}}
       \put(4,6.5){\line(-1,1){1}}\put(9,4.5){\line(1,0){2}}
       \put(9,4.5){\line(1,2){1}}\put(9,4.5){\line(0,1){2}}
       \put(9,6.5){\line(1,0){1}}\put(10,6.5){\line(1,-2){1}}
       \put(10,6.5){\line(1,-1){1}}\put(11,4.5){\line(0,1){1}}
       \put(1,1.6){$v_1$}\put(3.8,1.6){$v_2$}\put(8.4,2.2){$v_1$}\put(11.4,2.2){$v_2$}
       \put(7.6,6.4){$2$}\put(8.5,4){$k-1$}\put(9.9,4){$k$}\put(10.7,4){$k+1$}
       \put(2.4,1.87){$\bullet$}\put(9.93,-0.1){$\bullet$} \put(13,5.5){\rm
         (a)}\put(13,1.5){\rm
         (b)}\put(6.5,1.9){$\Longrightarrow$}\put(6,5.4){$\Longrightarrow$}
       \put(2.4,2.2){$\bx_0$}\put(10.2,0){$\widehat\bx_0$}
       \put(1.7,7.7){$\bm_c$}\put(4.1,6.4){$\bm_e$}
     \end{picture}
     \caption{Curves of type (b.1)}\label{fn2}
   \end{figure}
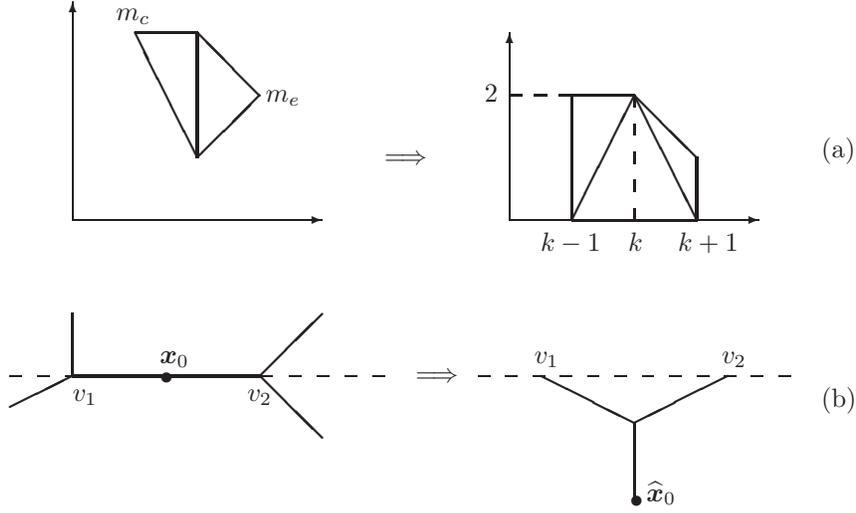

   Now we are going to refine the tropical limit of $C$ as described in
   \cite[Section 3.5]{Shu06a}. Geometrically it corresponds to a
   blow up resolving the considered singularity. Without loss of generality, assume that
   $\sigma$ is vertical (cf. Figure \ref{fig:differentsides}). Consider the polynomial
   \begin{equation}\widehat f(x,y)=f(x,y+1)\ .\label{e3}\end{equation}
   It defines a curve $\widehat
   C\subset\K\times\K^*$ with a singularity at $\widehat p=(1,0)$.
   Note that this refined tropical curve is a member of the family we
   described in Section \ref{sec-nottorus}, and indeed we will
   see that the refined tropical curve has a fat down end.
   Denote by $\widehat {\mathcal{C}}$ its tropicalisation and by
   $\widehat\nu:\widehat\Delta\to\R$ the corresponding concave piece-wise
   linear function on its Newton polygon $\widehat\Delta$. The fragment
   $(\delta',\delta'',C_{\delta'},C_{\delta''})$ of the tropical limit of $C$
   turns into the fragment
   $(\widehat\delta',\widehat\delta'',\delta_\sigma,\widehat C_{\delta'},\widehat
   C_{\delta''},C_\sigma)$ of the tropical limit of $\widehat C$ (see
   Figure \ref{fn2}(a)), where $C_\sigma$ is a curve in the toric surface
   $\Tor(\delta_\sigma)$, $\delta_\sigma=\conv\{(k-1,0),(k,2),(k+1,0)\}$,
   having a singularity at $\widehat p=(1,0)$. Observe that, since
   $\widehat p$ appears on the toric divisor $\Tor(\widehat\sigma)$,
   where $\widehat\sigma=[(k-1,0),(k+1,0)]$, and since $\widehat f(1,0)=0$,
   which implies $\ini_{(0,-\infty)}\widehat f(1,0)=0$ we can conclude that
   the values
   $\widehat\nu(k-1,0)$ and $\widehat\nu(k+1,0)$ must be equal.
   In
   view of the clear relations $\widehat\nu(k-1,0)=\nu(\bm_c)$ and
   $\widehat\nu(k+1,0)=\nu(\bm_e)$, this confirms the equality
   $\nu(\bm_c)=\nu(\bm_e)$, equivalent to the metric relation $l_1=l_2$
   for the tropical curve $\mathcal{C}$ shown in Section \ref{sec-maxclass}.

   Furthermore, from the above refinement we derive a correct canonical
   parametrisation of $\mathcal{C}$. It can be conveniently represented via the
   tropical blow-up\footnote{The above refinement can be interpreted as
     the (weighted) blow-up of the toric variety $\Tor(\widetilde\Delta)$
     at the point $p$ which replaces it by the exceptional divisor
     $\Tor(\delta_\sigma)$ (cf. \cite[Section 2]{ST06}).}, or {\it
     modification} in the terminology of \cite[Section 1]{Mi07a}. Notice
   that (under our assumptions) the edge $E$ of $\mathcal{C}$ dual to $\sigma$ lies
   on the $x$-axis of $\R^2$. Replace the (tropical) plane $\R^2$ by
   the tropical plane $P$ in $\R^3$ which consists of three
   half-planes:
   $$P_+=\{y\ge0,\ y=z\},\ P_-=\{y\le0,\ z=0\},\ P_0=\{y=0,\ z\le 0\}\
   ,$$ and introduce the projection $$\pi:P\to\R^2,\quad
   \pi(x,y,z)=(x,y)\ .$$ Then the tropical curve $\mathcal{C}\subset\R^2$ lifts
   to a tropical curve $\mathcal{C}^*\subset P$ 
   \begin{itemize}
   \item the part $\mathcal{C}\backslash E$ lifts to
     $\pi^{-1}(\mathcal{C}\backslash E)\cap(P_+\cup P_-)$,
   \item the edge $E$ is
     replaced by the fragment dual to the triangle $\delta_\sigma$ placed in
     $P_0$ (see Figure \ref{fn2}(b)), whereas the (tropical) singular
     point $\bx_0=(0,0)$ lifts to the (infinite) $1$-valent vertex
     $\widehat\bx_0$ of the vertical ray of the above
     fragment.
   \end{itemize}
   The map $\pi:\mathcal{C}^*\to \mathcal{C}$ provides a
   canonical parametrisation of $\mathcal{C}$ of genus $g$.

   \medskip

   \subsection{ Curves of type (b.2).} 
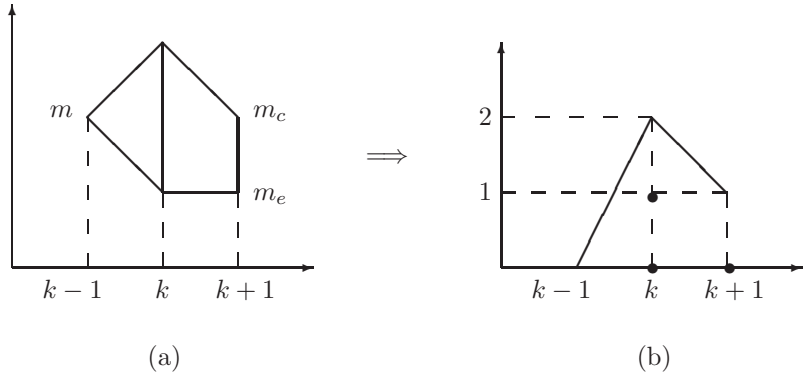
\begin{figure}[h]
     \setlength{\unitlength}{1cm}
     \begin{picture}(12,5)(0.3,0)
       \thinlines\put(1,1.5){\vector(0,1){3.5}}\put(1,1.5){\vector(1,0){4}}
       \put(7.5,1.5){\vector(0,1){3}}\put(7.5,1.5){\vector(1,0){4}}
       \dashline{0.2}(2,1.5)(2,3.5)\dashline{0.2}(3,1.5)(3,2.5)
       \dashline{0.2}(4,1.5)(4,2.5)\dashline{0.2}(7.5,3.5)(9.5,3.5)
       \dashline{0.2}(9.5,1.5)(9.5,3.5)\dashline{0.2}(10.5,1.5)(10.5,2.5)
       \dashline{0.2}(7.5,2.5)(10.5,2.5)
       \thicklines\put(2,3.5){\line(1,1){1}}\put(2,3.5){\line(1,-1){1}}
       \put(3,2.5){\line(0,1){2}}\put(3,2.5){\line(1,0){1}}
       \put(3,4.5){\line(1,-1){1}}\put(4,2.5){\line(0,1){1}}
       \put(8.5,1.5){\line(1,2){1}}\put(9.5,3.5){\line(1,-1){1}}
       \put(1.4,1.1){$k-1$}\put(2.9,1.1){$k$}\put(3.7,1.1){$k+1$}
       \put(7.2,3.4){$2$}\put(7.2,2.4){$1$}\put(7.9,1.1){$k-1$}
       \put(9.4,1.1){$k$}\put(10.2,1.1){$k+1$}
       \put(9.43,1.4){$\bullet$}\put(9.43,2.35){$\bullet$}\put(10.45,1.4){$\bullet$}
       \put(2.8,0.2){\rm (a)}\put(9.3,0.2){\rm
         (b)}\put(5.7,2.9){$\Longrightarrow$}\put(4.2,3.5){$\bm_c$}\put(4.2,2.4){$\bm_e$}
       \put(1.5,3.5){$m$}
     \end{picture}
     \caption{Curves of type (b.2), I}\label{fn3}
   \end{figure}

In this case, the edge $\sigma$ of the
   subdivision $S_\mathcal{C}$ spanned by the circuit $Z$ has lattice length $2$
   and is common for a triangle $\delta'$ and a trapeze $\delta''$,
   altogether containing $6$ integral points (see Figure \ref{fn3}(a)).
   Without loss of generality assume that $\sigma$ is vertical (see
   Figure \ref{fn3}(a)). Then the function $\nu:\Delta\to\R$ satisfies
   the following:
   \begin{equation}\begin{cases} &\nu\big|_\sigma=0,\quad \nu(k,j)<0,\ (k,j)\not\in\sigma,\\
       &\nu(m)=\alpha<0,\quad\nu(k-1,j)<\alpha,\ (k,j)\ne m,\\
       &\nu(\bm_c)=\nu(\bm_e)=\beta<0,\quad\nu(k+1,j)<\beta,\ (k,j)\ne \bm_c,\bm_e,\\
       &\nu(k+s,j)<s\beta,\quad s\ge 2\ .\end{cases}\label{e5}\end{equation}
   Due to the generality condition (G), we can assume that $\alpha<\beta$
   (the inequality $\alpha\le\beta$ is included in the definition of the
   corresponding cone of $\Trop(\Sing_\K(\Delta))$).

   The limit point $p$ belongs to $\Tor(\sigma)$, since the function
   $\nu:\Delta\to\R$ vanishes only along $\sigma$. Then the limit curve
   $C_{\delta'}$ is nonsingular and quadratically tangent to $\Tor(\sigma)$
   at $p$. The limit curve $C_{\delta''}$ either is
   nonsingular, quadratically tangent to $\Tor(\sigma)$ at $p$, or splits
   into two components transversally intersecting at $p$. The former
   option is not possible, since, otherwise, the refinement used in the
   preceding stage would lead to a subdivision containing the triangle
   $\delta_\sigma$ as in Figure \ref{fn2}(a), and hence to the equality
   $\alpha=\beta$ against the assumption made. Thus, $C_{\delta''}$ is
   reducible as indicated above.

   Let $\widehat f(x,y)=f(x,y+1)$. In view of (\ref{e5}), the fragment
   $(\delta',\delta'')$ of the subdivision $S_\mathcal{C}$ of $\Delta$ turns into a
   fragment of the subdivision $S_{\widehat {\mathcal{C}}}$ (in the notation of the
   preceding step) containing the two edges shown in Figure
   \ref{fn3}(b) with the following values of the function
   $\widehat\nu:\widehat\Delta\to\R$:
   \begin{equation}
     \widehat\nu(k-1,0)=\alpha,\
     \widehat\nu(k,0),\widehat\nu(k,1)<\widehat\nu(k,2)=0,\
     \widehat\nu(k+1,0)<\widehat\nu(k+1,1)=\beta\
     .\label{e4}
   \end{equation}
   To derive these relations, notice, first, that the expansion of
   $f(x,y)$ into power series in $t$ looks as
   $f(x,y)=x^ky^{k'}(y-1)^2+O(t)$, where $(k,k')$ is the bottom vertex
   of $\sigma$, and hence $\widehat f(x,y)=x^ky^2+O(y^3)+O(t)$. Second,
   recall that the limit curve $C_{\delta''}$ is reducible and both
   its components hit the point $(0,1)$, which is the limit point of
   $\bp=(1,1)$ in $\Tor(\sigma)$. The shape of
   $\delta''$ dictates that one of the components is $\{y-1=0\}$, in
   particular, the truncation of $f(x,y)$ to the edge $[m_c,m_e]$ is
   $x^{k+1}y^{k'}(y-1)t^\beta(1+O(t))$, and hence the substitution
   $y\to y+1$ produces the truncation $t^\beta x^{k+1}y(1+O(y)+O(t))$ of
   $\widehat f(x,y)$ to the segment $[(k+1,0),(k+1,1)]$.

   \medskip

\begin{figure}[h]
     \setlength{\unitlength}{0.85cm}
     \begin{picture}(13,4)(0,0)
       \thinlines\put(0.5,1.5){\vector(0,1){2.5}}\put(0.5,1.5){\vector(1,0){3.5}}
       \put(5,1.5){\vector(0,1){2.5}}\put(5,1.5){\vector(1,0){3.5}}
       \put(9.5,1.5){\vector(0,1){2.5}}\put(9.5,1.5){\vector(1,0){3.5}}
       \dashline{0.2}(0.5,3.5)(2.5,3.5)\dashline{0.2}(0.5,2.5)(3.5,2.5)
       \dashline{0.2}(2.5,1.5)(2.5,3.5)\dashline{0.2}(5,3.5)(7,3.5)
       \dashline{0.2}(5,2.5)(7,2.5)\dashline{0.2}(7,1.5)(7,2.5)
       \dashline{0.2}(9.5,3.5)(11.5,3.5)\dashline{0.2}(9.5,2.5)(12.5,2.5)
       \dashline{0.2}(11.5,1.5)(11.5,3.5)
       \thicklines\put(1.5,1.5){\line(1,2){1}}\put(2.5,3.5){\line(1,-1){1}}
       \put(3.5,1.5){\line(0,1){1}}\put(1.5,1.5){\line(1,0){2}}
       \put(6,1.5){\line(1,2){1}}\put(6,1.5){\line(1,1){1}}
       \put(6,1.5){\line(1,0){2}}\put(7,2.5){\line(1,0){1}}
       \put(7,2.5){\line(0,1){1}}\put(7,3.5){\line(1,-1){1}}
       \put(10.5,1.5){\line(1,0){2}}\put(10.5,1.5){\line(1,2){1}}
       \put(11.5,3.5){\line(1,-2){1}}\put(11.5,3.5){\line(1,-1){1}}
       \put(12.5,1.5){\line(0,1){1}}\put(8,1.5){\line(0,1){1}}
       \put(1.1,1.1){$k-1$}\put(2.4,1.1){$k$}\put(3.1,1.1){$k+1$}
       \put(5.6,1.1){$k-1$}\put(6.9,1.1){$k$}\put(7.6,1.1){$k+1$}
       \put(10.1,1.1){$k-1$}\put(11.4,1.1){$k$}\put(12.3,1.1){$k+1$}
       \put(0.2,3.4){$2$}\put(0.2,2.4){$1$}\put(4.7,3.4){$2$}\put(4.7,2.4){$1$}
       \put(9.2,3.4){$2$}\put(9.2,2.4){$1$}\put(2.3,0.2){\rm
         (a)}\put(6.8,0.2){\rm (b)}\put(11.3,0.2){\rm (c)}
     \end{picture}
     \caption{Curves of type (b.2), II}\label{fn4}
   \end{figure}
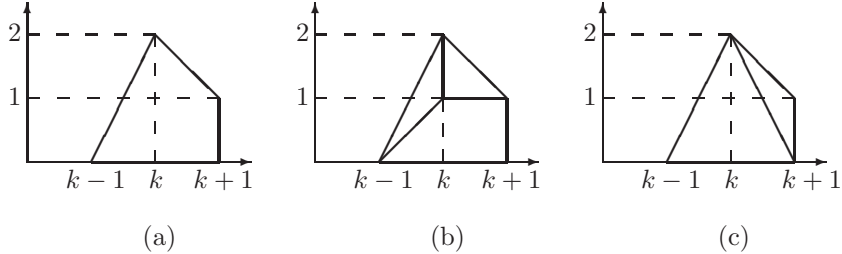

   {\it (1) Assume that the segment $[\bm_c,\bm_e]$ lies on
     $\partial\Delta$}. Then the subdivision $S_{\widehat {\mathcal{C}}}$ contains a
   fragment bounded by the quadrangle
   $$Q=\conv\{(k-1,0),(k,2),(k+1,1),(k+1,0)\}$$ (see Figure \ref{fn4}(a)).
   Since the point $\widehat \bp=(1,0)$ is singular for $\widehat C$,
   the limit point $\widehat p$ is singular for the corresponding
   limit curve of $\widehat C$. Thus,

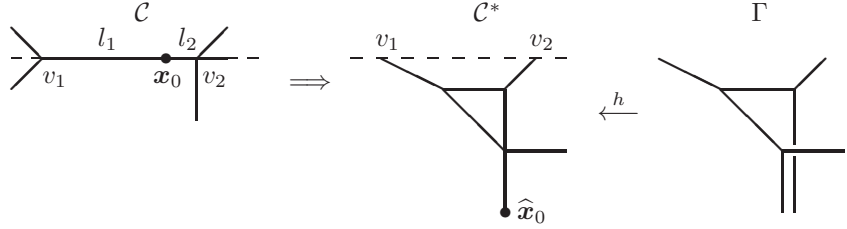
   \begin{figure}[h]
     \setlength{\unitlength}{0.82cm}
     \begin{picture}(13,10)(-0.2,0)
       \thinlines \dashline{0.2}(0,9)(0.5,9)\dashline{0.2}(3.5,9)(4,9)
       \dashline{0.2}(5.5,9)(9,9)\dashline{0.2}(0.5,3.5)(1,3.5)
       \dashline{0.2}(4.5,3.5)(5,3.5)\dashline{0.2}(7.5,3.5)(12.5,3.5)
       \thicklines\put(0,8.5){\line(1,1){0.5}}\put(0,9.5){\line(1,-1){0.5}}
       \put(0.5,9){\line(1,0){3}}\put(6,9){\line(2,-1){1}}
       \put(7,8.5){\line(1,-1){1}}\put(7,8.5){\line(1,0){1}}
       \put(8,6.5){\line(0,1){2}}\put(8,7.5){\line(1,0){1}}
       \put(8,8.5){\line(1,1){0.5}}\put(10.5,9){\line(2,-1){1}}
       \put(11.5,8.5){\line(1,0){1.2}}\put(11.5,8.5){\line(1,-1){1}}
       \put(12.5,6.5){\line(0,1){1}}\put(12.5,7.5){\line(1,0){1}}
       \put(12.7,6.5){\line(0,1){0.9}}\put(12.7,7.6){\line(0,1){0.9}}
       \put(12.7,8.5){\line(1,1){0.5}}\put(3,8){\line(0,1){1}}
       \put(3,9){\line(1,1){0.5}}\put(0.5,3){\line(1,1){0.5}}
       \put(0.5,4){\line(1,-1){0.5}}\put(1,3.5){\line(1,0){3.5}}
       \put(4,3.5){\line(1,1){0.5}}\put(4,3.5){\line(0,-1){1}}
       \put(8,3.5){\line(2,-1){2}}\put(10,1.5){\line(0,1){1}}
       \put(10,2.5){\line(2,1){1}}\put(11,3){\line(1,1){0.5}}
       \put(11,3){\line(1,0){1}} \put(0.5,8.6){$v_1$}\put(3.1,8.6){$v_2$}
       \put(2.4,8.9){$\bullet$}\put(2.3,8.6){$\bx_0$}\put(1.4,9.2){$l_1$}
       \put(2.7,9.2){$l_2$}\put(5.9,9.2){$v_1$}\put(8.4,9.2){$v_2$}
       \put(8.2,6.4){$\widehat\bx_0$}\put(7.9,6.4){$\bullet$}
       \put(2,9.6){$\mathcal{C}$}\put(7.5,9.6){$\mathcal{C}^*$}\put(12,9.6){$\Gamma$}
       \put(4.5,8.5){$\Longrightarrow$}\put(9.5,8){$\overset{h}{\longleftarrow}$}
       \put(1,3.1){$v_1$}\put(4.1,3.1){$v_2$}\put(2.9,3.38){$\bullet$}
       \put(2.8,3.1){$\bx_0$}\put(1.9,3.7){$l_1$}\put(3.4,3.7){$l_2$}
       \put(7.9,3.7){$v_1$}\put(11.4,3.7){$v_2$}\put(9.92,1.4){$\bullet$}
       \put(10.3,1.3){$\widehat\bx_0$}\put(2.5,4.1){$\mathcal{C}$}
       \put(9.5,4.1){$\mathcal{C}^*=\Gamma$}\put(6,3){$\Longrightarrow$}\put(4.5,5){$\text{(a):}\quad 
         l_2<l_1/2,\quad \alpha/2<\beta<0$}\put(4.5,0.5){$\text{(b):}\quad
         l_1/2<l_2\le l_1,\quad \alpha<\beta<\alpha/2$}
     \end{picture}
     \caption{Curves of type (b.2), III}\label{fn5}
   \end{figure}

   \begin{itemize}\item the entire segment $[(k-1,0),(k+1,0)]$ must be
     an edge of the induced subdivision of $Q$, and \item the limit curve
     $C'$ corresponding to the polygon which has the segment
     \mbox{$[(k-1,0),(k+1,0)]$} as a face must be singular at $\widehat p$.
   \end{itemize} This allows one only the following subdivisions of $Q$
   and relations on $\widehat\nu$: \begin{enumerate}\item[(i)] the
     subdivision shown in Figure \ref{fn4}(b), where the limit curve $C'$
     splits into two components transversally intersecting at $\widehat
     p$, and
     $$\widehat\nu(k-1,0)=\widehat\nu(k,0)=\widehat\nu(k+1,0)=\alpha,\
     \widehat\nu(k,1)=\widehat\nu(k+1,1)=\beta\in\left(\frac{\alpha}{2},0\right)\
     ;$$
   \item[(ii)] the subdivision shown in Figure \ref{fn4}(c), where the
     limit curve $C'$ is irreducible with a node at $\widehat p$, and
     $$\widehat\nu(k-1,0)=\widehat\nu(k,0)=\widehat\nu(k+1,0)=\alpha,\
     \widehat\nu(k,1)=\frac{\alpha}{2},\
     \widehat\nu(k+1,1)=\beta\in\left(\alpha,\frac{\alpha}{2}\right)\ .$$
   \end{enumerate} Notice that the relation $\alpha>\beta$ is not possible
   due to the last inequality in (\ref{e4}) which confirms the same
   conclusion of Section \ref{sec-maxclass}. We complete the study of this
   case with a canonical parametrisation of the tropical curve $\mathcal{C}$: we
   perform the modification of the plane as for the type (b.1) curves
   and replace the edge $E$ of $\mathcal{C}$ passing through the origin with the
   fragment dual to the subdivisions shown in Figure \ref{fn4}(b,c) -
   see Figure \ref{fn5}(a,b): in the first case, we have the
   parametrisation $\Gamma\overset{h}{\to}\mathcal{C}^*\overset{\pi}{\to}
   \mathcal{C}$, and, in the second case, the parametrisation $\Gamma=\widehat
{   \mathcal{C}}\overset{\pi}{\to} \mathcal{C}$. The geometry of those fragments of $\widehat
{   \mathcal{C}}$ imply the metric conditions on the position of the tropical
   singularity $\bx_0$ as indicated in Figure \ref{fn5}.

   \medskip

\begin{figure}[h]
     \setlength{\unitlength}{0.9cm}
     \begin{picture}(12,9)(0,0)
       \thinlines\put(1,6.5){\vector(0,1){2.5}}\put(1,6.5){\vector(1,0){5}}
       \put(1,1.5){\vector(0,1){2.5}}\put(1,1.5){\vector(1,0){5}}
       \put(7,6.5){\vector(0,1){2.5}}\put(7,6.5){\vector(1,0){5}}
       \put(7,1.5){\vector(0,1){2.5}}\put(7,1.5){\vector(1,0){5}}
       \dashline{0.2}(1,2.5)(4,2.5)\dashline{0.2}(1,3.5)(3,3.5)
       \dashline{0.2}(3,1.5)(3,3.5)\dashline{0.2}(1,7.5)(4,7.5)
       \dashline{0.2}(1,8.5)(3,8.5)\dashline{0.2}(3,6.5)(3,8.5)
       \dashline{0.2}(4,6.5)(4,7.5)\dashline{0.2}(7,2.5)(10,2.5)
       \dashline{0.2}(7,3.5)(9,3.5)\dashline{0.2}(10,1.5)(10,2.5)
       \dashline{0.2}(7,7.5)(10,7.5)\dashline{0.2}(7,8.5)(9,8.5)
       \dashline{0.2}(9,6.5)(9,7.5)
       \thicklines\put(2,6.5){\line(1,2){1}}\put(3,8.5){\line(1,-1){1}}
       \put(8,6.5){\line(1,2){1}}\put(8,6.5){\line(1,1){1}}
       \put(9,7.5){\line(0,1){1}}\put(9,8.5){\line(1,-1){2}}
       \put(9,7.5){\line(1,0){1}}\put(10,6.5){\line(0,1){1}}
       \put(2,1.5){\line(1,2){1}}\put(3,3.5){\line(1,-1){2}}
       \put(8,1.5){\line(1,2){1}}\put(8,1.5){\line(1,1){1}}
       \put(9,1.5){\line(0,1){2}}\put(9,2.5){\line(1,0){1}}
       \put(9,3.5){\line(1,-1){2}}\put(4,1.5){\line(0,1){1}}
       \put(3,3.5){\line(1,-2){1}}\put(8,6.5){\line(1,0){3}}
       \put(2,1.5){\line(1,0){3}}\put(8,1.5){\line(1,0){3}}
       \put(2.94,6.37){$\bullet$}\put(2.94,7.37){$\bullet$}
       \put(3.94,6.37){$\bullet$}\put(4.94,6.37){$\bullet$}
       \put(0.6,2.4){$1$}\put(0.6,3.4){$2$}\put(0.6,7.4){$1$}\put(0.6,8.4){$2$}
       \put(6.6,2.4){$1$}\put(6.6,3.4){$2$}\put(6.6,7.4){$1$}\put(6.6,8.4){$2$}
       \put(1.5,1.1){$k-1$}\put(2.9,1.1){$k$}\put(3.6,1.1){$k+1$}\put(4.8,1.1){$k+2$}
       \put(1.5,6.1){$k-1$}\put(2.9,6.1){$k$}\put(3.6,6.1){$k+1$}\put(4.8,6.1){$k+2$}
       \put(7.5,1.1){$k-1$}\put(8.9,1.1){$k$}\put(9.6,1.1){$k+1$}\put(10.8,1.1){$k+2$}
       \put(7.5,6.1){$k-1$}\put(8.9,6.1){$k$}\put(9.6,6.1){$k+1$}\put(10.8,6.1){$k+2$}
       \put(3,0.2){\rm (c)}\put(3,5.2){\rm (a)}\put(9,0.2){\rm
         (d)}\put(9,5.2){\rm (b)}
     \end{picture}
     \caption{Curves of type (b.2), IV}\label{fn6}
   \end{figure}
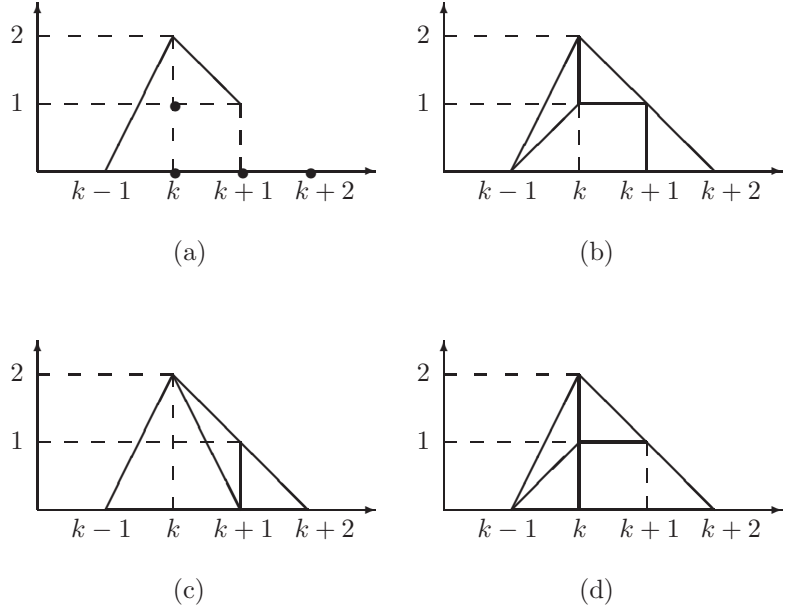

   {\it (2) Assume that the segment $[\bm_c,\bm_e]$ does not lie on
     $\partial\Delta$}. Then the fragment of the subdivision $S_{\widehat
{     \mathcal{C}}}$ we are interested in may include the point $(k+2,0)$ too (see
   Figure \ref{fn6}(a)), where due to (\ref{e5}),
   $\widehat\nu(k+2,0)<2\beta=2\widehat\nu(k+1,1)$. If
   $\widehat\nu(k+2,0)<\widehat\nu(k-1,0)=\alpha$, then the preceding
   argument leaves us with the only possible subdivisions shown in
   Figure \ref{fn6}(b,c) with the same conclusions as for the
   subdivisions in Figure \ref{fn4}(b,c) analysed before. If
   $\widehat\nu(k+2,0)>\alpha$ (what, in particular, yields
   $\beta>\alpha/2$), then the only suitable subdivision is shown in
   Figure \ref{fn6}(d), where
   $$\widehat\nu(k,0)=\widehat\nu(k+1,0)=\widehat\nu(k+2,0),\
   \widehat\nu(k,1)=\widehat\nu(k+1,1)=\beta\ ,$$ and the limit curve
   with the Newton trapeze splits into two components transversally
   intersecting the toric divisor $\Tor([(k,0),(k+2,0)])$ at the same
   point $\widehat\bp$.

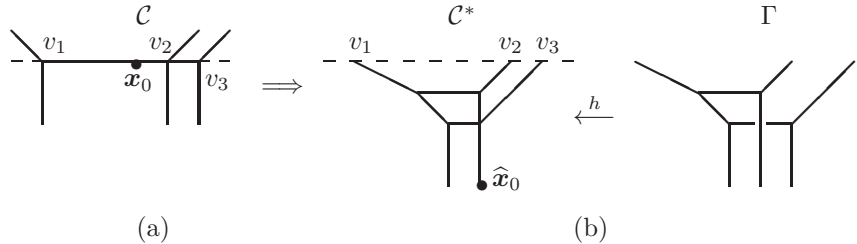
\begin{figure}[h]
     \setlength{\unitlength}{0.83cm}
     \begin{picture}(14,4)(-0.2,0)
       \thinlines \dashline{0.2}(0,3)(0.5,3)\dashline{0.2}(3,3)(3.5,3)
       \dashline{0.2}(5,3)(9,3)
       \thicklines\put(0.5,2){\line(0,1){1}}\put(0.5,3){\line(-1,1){0.5}}
       \put(0.5,3){\line(1,0){2.5}}\put(2.5,2){\line(0,1){1}}
       \put(3,2){\line(0,1){1}}\put(2.5,3){\line(1,1){0.5}}
       \put(3,3){\line(1,1){0.5}}\put(5.5,3){\line(2,-1){1}}
       \put(6.5,2.5){\line(1,-1){0.5}}\put(6.5,2.5){\line(1,0){1}}
       \put(7,1){\line(0,1){1}}\put(7,2){\line(1,0){0.5}}
       \put(7.5,1){\line(0,1){1.5}}\put(7.5,2.5){\line(1,1){0.5}}
       \put(7.5,2){\line(1,1){1}}\put(10,3){\line(2,-1){1}}
       \put(11,2.5){\line(1,0){1}}\put(11,2.5){\line(1,-1){0.5}}
       \put(11.5,1){\line(0,1){1}}\put(11.5,2){\line(1,0){0.4}}
       \put(12.1,2){\line(1,0){0.4}}\put(12,2.5){\line(1,1){0.5}}
       \put(12.5,2){\line(1,1){1}}\put(12,1){\line(0,1){1.5}}\put(12.5,1){\line(0,1){1}}
       \put(1.9,2.85){$\bullet$}\put(7.43,0.9){$\bullet$}\put(0.5,3.2){$v_1$}
       \put(2.2,3.2){$v_2$}\put(1.8,2.6){$\bx_0$}\put(5.4,3.2){$v_1$}\put(3.1,2.6){$v_3$}
       \put(7.8,3.2){$v_2$}\put(8.4,3.2){$v_3$}\put(7.7,1){$\widehat\bx_0$}
       \put(2,3.6){$\mathcal{C}$}\put(7,3.6){$\mathcal{C}^*$}\put(12,3.6){$\Gamma$}
       \put(4,2.5){$\Longrightarrow$}\put(9,2){$\overset{h}{\longleftarrow}$}
       \put(2,0.2){\rm (a)}\put(9,0.2){\rm (b)}
     \end{picture}
     \caption{Curves of type (b.2), V}\label{fn7}
   \end{figure}

   The canonical parametrisation of $\mathcal{C}$ again is built in the form
   $\Gamma\overset{h}{\to}\mathcal{C}^*\overset{\pi}{\to}\mathcal{C}$, where $\widehat
{   \mathcal{C}}$ appears in the modification of the plane along the $x$-axis: the
   edge $E$ of $\mathcal{C}$ (with the endpoints $v_1,v_3$ in Figure
   \ref{fn7}(a)) is replaced by the fragment dual to the subdivision in
   Figure \ref{fn6}(d) and lying in the half-plane $P_0$ (shown in
   Figure \ref{fn7}(b)). At last, the parameterising graph $\Gamma$ is
   obtained by the resolving the double ray with the endpoint
   $\widehat\bx$ (see Figure \ref{fn7}(c)). The singular point lifts to two $1$-valent vertices of the parametrising graph. Notice that this case
   corresponds to the metric relation $l_2<l_1/2$.

\begin{thebibliography}{10}

\bibitem{AK06}
Federico Ardila and Carly Klivans.
\newblock The {Bergman} complex of a matroid and phylogenetic trees.
\newblock {\em J. Combin. Theory Ser. B}, 96:38--49, 2006.

\bibitem{BJSST07}
Tristram Bogart, Anders Jensen, David Speyer, Bernd Sturmfels, and Rekha
  Thomas.
\newblock Computing tropical varieties.
\newblock {\em J. Symbolic Comput.}, 42:54--73, 2007.
\newblock arXiv:math.AG/0507563.

\bibitem{DFS05}
Alicia Dickenstein, Eva~Maria Feichtner, and Bernd Sturmfels.
\newblock Tropical discriminants.
\newblock {\em J. Amer. Math. Soc.}, 20:1111--1133, 2007.
\newblock arXiv:math.AG/0510126.

\bibitem{EKL06}
Manfred Einsiedler, Mikhail Kapranov, and Douglas Lind.
\newblock Non-archimedean amoebas and tropical varieties.
\newblock {\em J.\ Reine Angew.\ Math.}, 601:139--157, 2006.

\bibitem{FS05}
Eva~Maria Feichtner and Bernd Sturmfels.
\newblock Matroid polytopes, nested sets and {Bergman} fans.
\newblock {\em Portugaliae Mathematica}, 62:437--468, 2005.
\newblock arXiv:math.CO:0411260.

\bibitem{GKZ}
Israel~M. Gelfand, Mikhail~M. Kapranov, and Andrei~V. Zelevinsky.
\newblock {\em Discriminants, resultants, and multidimensional determinants}.
\newblock Birkh\"auser Boston Inc., 1994.

\bibitem{MS09}
Diane Maclagan and Bernd Sturmfels.
\newblock Introduction to tropical geometry.
\newblock book in preparation, available at
  http:\\www.warwick.ac.uk/staff/D.Maclagan/papers/TropicalBook.pdf.

\bibitem{Mi03}
Grigory Mikhalkin.
\newblock Enumerative tropical geometry in {${\mathbb{R}^2}$}.
\newblock {\em J. Amer. Math. Soc.}, 18:313--377, 2005.
\newblock arXiv:math.AG/0312530.

\bibitem{Mi07a}
Grigory Mikhalkin.
\newblock Introduction to tropical geometry.
\newblock Notes from the IMPA Lectures in Summer 2007, arXiv:0709.1049, 2007.

\bibitem{RST03}
J\"urgen Richter-Gebert, Bernd Sturmfels, and Thorsten Theobald.
\newblock First steps in tropical geometry.
\newblock {\em Idempotent Mathematics and Mathematical Physics, Proceedings
  Vienna}, 2003.
\newblock arXiv:math/0306366.

\bibitem{Shu06a}
Eugenii Shustin.
\newblock A tropical approach to enumerative geometry.
\newblock {\em St.\ Petersburg Math.\ J.}, 17:343--375, 2006.

\bibitem{ST06}
Eugenii Shustin and Ilya Tyomkin.
\newblock Patchworking singular algebraic curves, ii.
\newblock {\em Israel J.\ Math.}, 151:145--166, 2006.

\bibitem{Spe05}
David Speyer.
\newblock {\em Tropical Geometry}.
\newblock PhD thesis, University of California, Berkeley, 2005.

\bibitem{Stu02}
Bernd Sturmfels.
\newblock {\em Solving Systems of linear equations}.
\newblock CBMS Regional Conference Series in Mathematics 97, American
  Mathematical Society, Providence, 2002.

\bibitem{Tyo09}
Ilya Tyomkin.
\newblock Tropical geometry and correspondence theorems via toric stacks.
\newblock Preprint, 2009.

\end{thebibliography}
\end{document}